\documentclass[12pt]{amsart}

\usepackage{amssymb}
\usepackage{amscd}

\newcommand{\bt}{\begin{theorem}}
\newcommand{\et}{\end{theorem}}
\newcommand{\bp}{\begin{proposition}}
\newcommand{\ep}{\end{proposition}}

\newcommand{\bl}{\begin{lemma}}
\newcommand{\el}{\end{lemma}}

\newcommand{\br}{\begin{result}}
\newcommand{\er}{\end{result}}
\newcommand{\be}{\begin{equation}}
\newcommand{\ee}{\end{equation}}
\newcommand{\bc}{\begin{corollary}}
\newcommand{\ec}{\end{corollary}}
\newcommand{\bex}{\begin{example}}
\newcommand{\eex}{\end{example}}

\newtheorem{theorem}{Theorem}[section]
\newtheorem{corollary}[theorem]{Corollary}
\newtheorem{lemma}[theorem]{Lemma}
\newtheorem{proposition}[theorem]{Proposition}
\newtheorem{result}[theorem]{Result}
\newtheorem{example}[theorem]{Example}

\newcommand{\K}{\mathbb{K}}
\newcommand{\N}{\mathbb{N}}

\newcommand{\cPA}{\mathcal{P}\!\mathcal{A}}
\newcommand{\cPAi}{\mathcal{P}\!\mathcal{A}\mathit{i}}

\newcommand{\cB}{\mathcal{B}}

\newcommand{\cD}{\mathcal{D}}
\newcommand{\cH}{\mathcal{H}}
\newcommand{\cI}{\mathcal{I}}
\newcommand{\cJ}{\mathcal{J}}
\newcommand{\cK}{\mathcal{K}}
\newcommand{\cL}{\mathcal{L}}
\newcommand{\cR}{\mathcal{R}}

\newcommand{\cSb}{\mathcal{S}\!\mathit{ub}}
\newcommand{\cSbi}{\mathcal{S}\!\mathit{ubi}\,}

\newcommand{\tsq}{\,${\tiny$\square$}$\,}

\newcommand{\lis}{[\![}
\newcommand{\ris}{]\!]}

\newcommand{\rg}{{\rm Reg}\,}

\newcommand{\Ng}{{\rm N}_}

\newcommand{\m}{{\rm M}_}

\newcommand{\med}{\medskip}
\newcommand{\sm}{\smallskip}

\begin{document}

\noindent {\em Algebra Universalis} {\bf 53} (2005), 407 -- 432\\
DOI 10.1007/s00012-005-1910-8
\vspace{0.05in}\\
{\em arXiv version}: layout, fonts, pagination and numbering of sections and theorems may vary from the AU published version
\vspace{0.1in}\\

\title{$\cPA$-isomorphisms of inverse semigroups}
\author{Simon M. Goberstein}
\address{\hspace*{-0.205in}Department of Mathematics \!\&\! Statistics,
California State University, Chico, CA 95929, USA, \!\!e-mail: sgoberstein@csuchico.edu}
\begin{abstract}
A partial automorphism of a semigroup $S$ is any isomorphism between its subsemigroups, and the set all partial automorphisms of
$S$ with respect to composition is the inverse monoid called the partial automorphism monoid of $S$. Two semigroups are said to be $\cPA$-isomorphic if their partial automorphism monoids are isomorphic. A class $\K$ of semigroups is  called $\cPA$-closed if it contains every semigroup $\cPA$-isomorphic to some semigroup from $\K$. Although the class of all inverse
semigroups is not $\cPA$-closed, we prove that the class of inverse semigroups, in which no maximal isolated subgroup is a
direct product of an involution-free periodic group and the two-element cyclic group, is $\cPA$-closed. It follows 
that the class of all combinatorial inverse semigroups (those with no nontrivial subgroups) is $\cPA$-closed. A semigroup is
called $\cPA$-determined if it is isomorphic or anti-isomorphic to any semigroup that is $\cPA$-isomorphic to it. We show
that combinatorial inverse semigroups which are either shortly connected [5] or quasi-archimedean [10] are $\cPA$-determined.\\

\noindent 2000 Mathematics Subject Classification: 20M10, 20M18, 20M20
\end{abstract}

\maketitle

\font\caps=cmcsc10 scaled \magstep1   

\setcounter{section}{-1}

\section{Introduction}
\sm Let $A$ be an algebraic structure of a certain type (e.g., a ring, a group, a semigroup, etc.); we will call it, for short,
an algebra and refer to its substructures as subalgebras of $A$. In many cases, it is convenient to regard the empty set as a
subalgebra of $A$ (especially if the intersection of two nonempty subalgebras of $A$ may be empty), and we will adhere to this
convention if $A$ is a semigroup or an inverse semigroup. A partial automorphism of $A$ is any isomorphism between its
subalgebras, and the set of all partial automorphisms of $A$ with respect to composition is an inverse monoid called the partial
automorphism monoid of $A$. The problem of characterizing algebras of various types by their partial automorphism monoids was
posed by Preston in [17]. It has been considered in a number of publications for several classes of groups and semigroups. In
[5] the present author studied  the problem of characterizing inverse semigroups $S$ by their partial automorphism monoids,
composed of all isomorphisms between {\em inverse subsemigroups} of $S$, in the class of {\em all inverse semigroups}. In this
article we investigate to what extent an inverse semigroup $S$ is characterized by its partial automorphism monoid, consisting
of all isomorphisms between {\em subsemigroups} of $S$, in the class of {\em all semigroups}.

The main results of the paper are contained in Sections 3 and 4. Since the two-element cyclic group is clearly $\cPA$-isomorphic
to the two-element null semigroup, the classes of groups and of inverse semigroups are not $\cPA$-closed. We show that, in some
sense, this is the only ``anomaly'' -- according to Theorem 3.13 of Section 3, the class of all those inverse semigroups, in which
no maximal isolated subgroup is a direct product of $C_2$ and a periodic group containing no elements of even order, is
$\cPA$-closed. It follows that several other large classes of inverse semigroups are $\cPA$-closed (Corollary 3.14), including
the class of all combinatorial inverse semigroups. Combining this fact with some earlier results from [5] and [10], we prove in
Section 4 (Theorem 4.7) that combinatorial inverse semigroups are $\cPA$-determined if they are either shortly connected [5] or
faintly archimedean [10] (for definitions see Section 4). We also show (Examples 4.5 and 4.6) that there exist shortly connected
combinatorial inverse semigroups which are not faintly archimedean.

We use [3] and [7] as standard references for the algebraic theory of semigroups and, in general, follow the terminology and
notation of these monographs. For an extensive treatment of the theory of inverse semigroups we refer to [16]. However, for the
reader's convenience, the basic semigroup-theoretic concepts and facts used in the paper are reviewed in Section 1 and all the
necessary preliminaries on lattice isomorphisms and $\cPA$-isomorphisms of semigroups are included in Section 2. This makes the
paper essentially self-contained.

The main results of this paper were reported at the International Conference on Algebra in Honor of Ralph McKenzie held at
Vanderbilt University on May 21-24, 2002. \med
\def\bfseries{\normalsize\caps}
\section{The background}
\med Let $S$ be an arbitrary semigroup. An element $x\in S$ is called {\em regular} (in the sense of von Neumann for rings) if
there is $y\in S$ such that $xyx=x$. If $x, y\in S$ satisfy $xyx=x$, then for $x'=yxy$ we have $xx'x=x$ and $x'xx'=x'$, in which
case $x'$ is called an {\em inverse} of $x$. Thus an element of $S$ is regular if and only if it has an inverse in $S$ (in general, more than one). Denote by $\rg (S)$ the set of all regular elements of $S$. For any $A\subseteq S$, the set of all those idempotents of $S$ which are contained in $A$ will be denoted by $E_A$. In particular, $E_S$ is the set of all idempotents of $S$. It is clear that $E_S\not=\emptyset$ if and only if $\rg (S)\not=\emptyset$. A semigroup $S$ is called {\em regular} if $\rg (S)=S$, and {\em idempotent-commutative} if $E_S\not=\emptyset$ and $ef=fe$ for all $e,f\in E_S$. It is easily shown [22, Theorem 3.1] that every regular element $x$ of an idempotent-commutative semigroup has a unique inverse (denoted usually by $x^{-1}$).

A regular idempotent-commutative semigroup is called an {\em inverse semigroup}, and an {\em inverse monoid} is an inverse semigroup with an identity element. Thus every element of an inverse semigroup has exactly one inverse. Note that if $S$ is an idempotent-commutative semigroup, it is immediate from [22, Theorem 3.2] that $\rg (S)$ is the largest inverse subsemigroup of $S$. If $S$ is an inverse semigroup, set $x\leq y$ if and only if $x=xx^{-1}y$ for $x, y\in S$; then $\leq$ is a partial order relation on $S$, compatible with the operations of multiplication and inversion, which is called the {\em natural} {\em order relation} on $S$. Clearly, if $S$ is an inverse semigroup, $E_S$ is a semilattice whose natural order relation is the restriction to $E_S$ of the natural order relation on $S$. Let $E$ be a semilattice. In what follows, we will have an occasion to consider simultaneously $(E,\leq)$ and $(E,\leq^d)$ where $\leq^d$ denotes the partial order on $E$ dual to $\leq$; to shorten notation, we will write $E$ instead of $(E,\leq)$ and $E^d$ instead of $(E,\leq^d)$.

Let $X$ be any set. The {\em symmetric inverse semigroup} $\cI_X$ on $X$ is the inverse monoid under composition consisting of all {\em partial bijections} of $X$ (that is, all bijections between various subsets of $X$, including $\emptyset$). It is easily seen that the natural order relation on $\cI_X$ is precisely the extension $\subseteq$ of partial bijections of $X$ and that the idempotents of $\cI_X$ are the identity mappings $1_A\!: a\mapsto a\;(a\in A)$ where $A$ is an arbitrary subset of $X$ (we will not distinguish $1_A$ from the identity relation $\{(a, a)\,\vert\,a\in A\}$ on $A$). Note that $1_\emptyset=\emptyset$ and $1_A\circ 1_B=1_{A\cap B}$ for any $A,\,B\subseteq X$. It follows that the semilattice of idempotents of $\cI_X$ is actually a lattice isomorphic to the lattice of all subsets of $X$. By the Wagner-Preston representation theorem [16, Theorem IV.1.6], each inverse semigroup $S$ is isomorphically embeddable into $\cI_S$ and the natural order relation $\leq$ on $S$ corresponds to $\subseteq$ under this embedding. 

Let $X, Y, X'$, and $Y'$ be any sets, $\rho\subseteq X\times Y$, and $\rho'\subseteq X'\times Y'$. As in [23], we define 
$\rho\tsq\rho'\subseteq (X\times X')\times (Y\times Y')$ as follows: $((x, x'), (y, y'))\in\rho\tsq\rho'$ if and only if $(x, y)\in\rho$ and $(x', y')\in\rho'$. We will often encounter the situation when $X=X'$, $Y=Y'$, and $\varphi$ is a certain bijection of $X$ onto $Y$. In this case, it is clear that $\varphi\tsq\varphi$ is a bijection of $\cI_X$ onto $\cI_Y$, and $\alpha(\varphi\tsq\varphi)=\varphi^{-1} \circ\alpha\circ\varphi$ for any $\alpha\in\cI_X$.

Let $S$ be a semigroup. For $a,b\in S$, let $a\cL b\;[a\cR b,\,a\cJ b]$ if and only if $a$ and $b$ generate the same principal left [right, two-sided] ideal of $S$. Set $\cH=\cL\cap\cR$ and $\cD=\cL\vee\cR$. Thus $\cH\subseteq\cL\subseteq\cD,\;\cH\subseteq\cR\subseteq\cD$ and $\cD\subseteq\cJ$. The equivalences $\cH,\cL,\cR,\cD$ and $\cJ$ are called the {\em Green's relations} on $S$ [3, Chapter 2]. It is easily seen that $S$ is regular [inverse] if and only if each $\cL$-class and each $\cR$-class of $S$ contains at least one [exactly one] idempotent. For $\cK\in\{\cH,\cL,\cR,\cD,\cJ\}$, denote by $K_x$ the $\cK$-class of $S$ containing $x\in S$. Note that if $x\in\rg (S)$, then every element of the $\cD$-class $D_x$ is regular [3, Theorem 2.11 (i)]. Thus if $D$ is a $\cD$-class of $S$, then either no element of $D$ is regular or all elements of $D$ are regular; in the latter case, we say that $D$ is a {\em regular $\cD$-class} of $S$. If $D$ is a regular $\cD$-class of $S$, each
$\cL$-class and each $\cR$-class in $D$ contains at least one idempotent [3, Theorem 2.11 (ii)], and if, in addition, it is assumed that $S$ is an idempotent-commutative semigroup, it is clear that each $\cL$-class and each $\cR$-class in $D$ contains exactly one idempotent.

Let $S$ be a semigroup and $U$ a subsemigroup of $S$. We will use the superscript $U$ for the Green's relations on $U$ in order to distinguish them from the corresponding relations on $S$ (which we will write without superscripts). If
$\cK\in\{\cH,\cL,\cR,\cD,\cJ\}$, it is clear that $\cK^U\subseteq\cK\cap(U\times U)$. In general, these inclusions may be proper for every $\cK\in\{\cH,\cL,\cR,\cD,\cJ\}$. However, if $\rg (S)$ is a subsemigroup of $S$ (in particular, if $S$ is idempotent-commutative), it is immediate that if $U=\rg (S)$, $a\in U$, and $\cK\in\{\cH, \cL, \cR,\cD\}$, then $K^U_a=K_a$.

Let $S$ be a semigroup. Denote by $J(x)$ the principal ideal of $S$ generated by $x\in S$. The set of $\cJ$-classes of $S$
is partially ordered by the relation $\leq$ defined as follows: $J_x\leq J_y$ if and only if $J(x)\subseteq J(y)$ for $x, y\in
S$. Similarly one can partially order the set of $\cL$-classes and the set of $\cR$-classes of $S$. We say that $x\in S$ is a {\em
group element} of $S$ if it belongs to some subgroup of $S$; otherwise $x$ is a {\em nongroup element} of $S$. Thus $x\in S$ is
a group element if and only if $x\in H_e$ for some $e\in E_S$, and a nongroup element if and only if either $x\not\in\rg (S)$ or
$x\in\rg (S)$ but $xx'\not=x'x$ where $x'$ is some (any) inverse of $x$ in $S$. If $A\subseteq S$, denote by $\Ng A$ the set of
all nongroup elements of $S$ contained in $A$ (so, in particular, $\Ng S$ is the set of all nongroup elements of $S$).
Let $D$ be a $\cD$-class of $S$. We say that $D$ is a {\em nongroup $\cD$-class} if $\Ng D\not=\emptyset$; otherwise $D$ is a
{\em group $\cD$-class}. Following Jones, we will also say that an idempotent $e$ of $S$ (and each subgroup of $H_e$) is {\em
isolated} if $D_e=H_e$, and {\em nonisolated} otherwise (see, for example, [9, p. 325] where these terms were introduced for
inverse semigroups). Thus an idempotent $e$ of $S$ is isolated if and only if $D_e$ is a group $\cD$-class.

Let $T$ be any semigroup with zero and $T^*=T\setminus\{0\}$. Take an arbitrary semigroup $A$ disjoint from $T^*$, and let $\eta : T^*\rightarrow A$ be a partial homomorphism (that is, $(xy)\eta=(x\eta)(y\eta)$ whenever $x, y, xy\in T^*$).  Denote $S=A\cup T^*$. For any $x, y\in S$, define $x\circ y$ as follows: $x\circ y=x(y\eta)$ if $x\in A$, $y\in T^*$; $x\circ y=(x\eta)y$ if $x\in T^*$, $y\in A$; $x\circ y=(x\eta)(y\eta)$ if $x, y\in T^*$ and $xy=0$ in $T$; and if $x, y\in T^*$ and $xy\in T^*$ or if $x, y\in A$, then $x\circ y$ coincides with the product of $x$ and $y$ in $T$ or $A$, respectively. Then $(S, \circ)$ is a semigroup whose operation is {\em determined by the partial homomorphism $\eta$}, and $S$ is a {\em retract ideal extension} of $A$ by $T$ [16, \S I.9]; conversely, if a semigroup is a retract ideal extension of $A$ by $T$, then its operation is determined by some partial homomorphism of $T^*$ to $A$ [16, Proposition I.9.14]. In what follows the word ``extension'' will be used instead of ``retract ideal extension'' since we will be dealing only with such ideal extensions of semigroups.

Let $S$ be an arbitrary semigroup. If $X$ is a nonempty subset of $ S$, the subsemigroup of $S$ generated by $X$ will be denoted by $\langle X\rangle$. Take any $x\in S$. Then $\langle x\rangle$ is called the {\em monogenic subsemigroup
of $S$ generated by $x$}. If $\langle x\rangle$ is finite, the {\em order} of $x$ is the number of elements of $\langle
x\rangle$; it will be denoted by $o(x)$. If $\langle x\rangle$ is infinite, we write $o(x)=\infty$ and say that $x$ has {\em
infinite order}. If $o(x)<\infty$, the {\em index} of $x$ (to be denoted by $\text{ind }x$) is defined as the least positive
integer $m$ satisfying $x^m=x^{m+k}$ for some positive integer $k$, and the smallest of such integers $k$ is called the {\em
period} of $x$ [7, p. 8]. If $x$ has infinite order, we set $\text{ind }x=\infty$. If $x$ has finite order and if $m$ and $n$
stand for its index and period, respectively, the monogenic semigroup $\langle x\rangle$ (or any semigroup isomorphic to it)
will be denoted by $M(m, n)$. It is easily seen that for any $m, n\in\N$, there is one and (up to isomorphism) only one
monogenic semigroup $M(m, n)$ [7, Section I.2]. Clearly, $M(1, n)$ is the cyclic group of order $n$ for which we will adopt the
commonly used notation $C_n$. Denote by $\m S$ the set of all $x\in S$ such that the monogenic semigroup $\langle x\rangle$ has
a unique generator. It is immediate that $x\in\m S$ if and only if either $\text{ind }x=1$ and $o(x)\leq 2$ or $\text{ind }x>1$.
Thus $\Ng S\cup E_S\subseteq \m S$, and if $x$ is a nonidempotent group element of $\m S$, then  $o(x)=2$ or $o(x)=\infty$.
Recall also that a semigroup $\langle a,b\rangle$ with identity $1$ given by one defining relation $ab=1$ is said to be {\em
bicyclic} [3, \S\,1.12]; we will denote it by $\cB(a, b)$. The idempotents of $\cB(a, b)$ form a chain: $1=ab>ba>b^2a^2>\ldots$,
and $\cB(a, b)$ is an inverse monoid consisting of a single $\cD$-class [3, Theorem 2.53]. A semigroup $S$ is said to be {\em
completely semisimple} if it contains no bicyclic subsemigroup.

Let $S$ be an inverse semigroup. If $X$ is a nonempty subset of $S$, the {\em inverse} subsemigroup of $S$ generated by $X$ will
be denoted by $\lis X \ris$, so $\lis X \ris=\langle X\cup X^{-1}\rangle=\lis X^{-1}\ris$ where $X^{-1}=\{x^{-1}\,\vert\,x\in
X\}$. If $X=\{x\}$ for some $x\in S$, we will write $\lis x\ris$ instead of $\lis X\ris$ and call $\lis x\ris$ the {\em
monogenic inverse subsemigroup of $S$ generated by $x$}; if $S=\lis x\ris$, we will say that the inverse semigroup $S$ is {\em
monogenic}. A detailed analysis of the structure of monogenic inverse semigroups is contained in [16, Chapter IX]). We recall
only a few basic facts about them. Let $S=\lis x\ris$ be a monogenic inverse semigroup. Then $\cD=\cJ$ and the partially ordered
set of $\cJ$-classes (=$\cD$-classes) of $S$ is a chain with the largest element $J_x\;(=D_x)$. It is obvious that one of the following
holds: (a) $xx^{-1}=x^{-1}x$, (b) $xx^{-1}\,\|\, x^{-1}x$, (c) $xx^{-1}>x^{-1}x$ or $x^{-1}x>xx^{-1}$. In case (a), $S=D_x$ is a
cyclic group. In case (b), $D_x=\{x,x^{-1},xx^{-1},x^{-1}x\}\;(=J_x)$ is the ``top'' $\cJ$-class of $S$ (and, of course, $S\setminus D_x$
is an ideal of $S$). Finally, in case (c), $S=D_x$ is a bicyclic semigroup -- it is either $\cB(x,x^{-1})$ if $xx^{-1}>x^{-1}x$,
or $\cB(x^{-1},x)$ if $x^{-1}x>xx^{-1}$. \med
\section{Lattice isomorphisms and $\cPA$-isomorphisms of semigroups, preliminaries}
\med Let $S$ be a semigroup. Since we assume that $\emptyset$ is a subsemigroup of $S$, the set of all subsemigroups of $S$,
partially ordered by inclusion, is a complete lattice which we will denote by $\cSb(S)$. It is clear that $H\cap K$ is the
greatest lower bound and $\langle H\cup K\rangle$ is the least upper bound of $H,\,K\in\cSb(S)$; we will usually denote the
latter by $H\vee K$. Let $T$ be a semigroup such that $\cSb(S)\cong\cSb(T)$. Then $S$ and $T$ are said to be {\em lattice
isomorphic}, and any isomorphism of $\cSb(S)$ onto $\cSb(T)$ is called a {\em lattice isomorphism} of $S$ onto $T$.  If $\Psi$
is a lattice isomorphism of $S$ onto $T$, we say that $\Psi$ is {\em induced} by a mapping $\psi\!: S\rightarrow T$ (or that
$\psi$ {\em induces} $\Psi$) if $H\Psi=H\psi$ for all $H\in\cSb(S)$.

Let $S$ and $T$ be lattice isomorphic semigroups, and let $\Psi$ be an isomorphism of $\cSb(S)$ onto $\cSb(T)$. It is obvious
(and well-known [21, Lemma 3.1(b)]) that a subsemigroup $U$ of $S$ is an atom of $\cSb(S)$ if and only if $U=\langle
e\rangle=\{e\}$ for some idempotent $e\in S$. Thus $E_S\not=\emptyset$ if and only if $E_T\not=\emptyset$, and there is a unique
bijection $\psi_E$ of $E_S$ onto $E_T$ defined by the formula $\{e\}\Psi=\{e\psi_E\}$ for all $e\in E_S$. We will say that
$\psi_E$ is the {\em $E$-bijection associated with $\Psi$}. It is also easily seen (and well-known [21, Proposition 36.6]) that
for all $e, f\in E_S$, we have $e\nparallel f$ if and only if $e\psi_E\nparallel f\psi_E$, and if $e\|f$, then
$(ef)\psi_E=(e\psi_E)(f\psi_E)$, which is expressed by saying that $\psi_E$ is a {\em weak isomorphism} of $E_S$ onto $E_T$.

A {\em partial automorphism} of a semigroup $S$ is any isomorphism between its subsemigroups. We denote by $\cPA(S)$ the set of
all partial automorphisms of $S$. Since $\emptyset\in\cSb(S)$, it is natural to regard $\emptyset$ as the (unique) isomorphism
of the empty subsemigroup of $S$ onto itself, so $\emptyset\in\cPA(S)$. With respect to composition $\cPA(S)$ is an inverse semigroup which is an inverse subsemigroup of $\cI_S$. In particular, the natural order relation on $\cPA(S)$ coincides with the extension 
$\subseteq$ of partial automorphisms of $S$ and the idempotents of $\cPA(S)$ are precisely the identity mappings 
$1_H\!: h\mapsto h\;(h\in H)$ where $H\in\cSb(S)$. Clearly, $1_\emptyset\;(=\emptyset)$ is the zero while $1_S$ is the identity of $\cPA(S)$. Thus $\cPA(S)$ is an inverse monoid with zero; it is called the {\em partial automorphism monoid} of $S$. The group of units of $\cPA(S)$ is the automorphism group of $S$, and the semilattice of idempotents of $\cPA(S)$ is a lattice isomorphic to $\cSb(S)$.

Let $S$ and $T$ be semigroups. If $\cPA(S)\cong\cPA(T)$, then $S$ and $T$ are said to be $\cPA${\em-isomorphic}, and any
isomorphism of $\cPA(S)$ onto $\cPA(T)$ is called a $\cPA${\em-isomorphism} of $S$ onto $T$. Let $\Phi$ be a $\cPA$-isomorphism
of $S$ onto $T$. For any $H\in\cSb(S)$, define $H\Phi^*$ by the formula $1_H\Phi=1_{H\Phi^*}$. Then $\Phi^*$ is a lattice
isomorphism of $S$ onto $T$. If $E_S\not=\emptyset$, we will denote by $\varphi_E$ the $E$-bijection associated with $\Phi^*$
and say that it is {\em associated with $\Phi$}; thus $\{e\}\Phi^*=\{e\varphi_E\}$ for all $e\in E_S$. If there is a bijection
$\varphi\!:S\rightarrow T$ such that $\Phi=(\varphi\tsq\varphi)\vert_{\cPA(S)}$, we say that $\varphi$ {\em induces} $\Phi$ (or
$\Phi$ is {\em induced} by $\varphi$). Thus $\Phi$ is induced by $\varphi$ if for all $\alpha\in\cPA(S)$ and $x,y\in S$, we have
$\;x\alpha=y$ if and only if $(x\varphi)(\alpha\Phi)=y\varphi$. Let $\theta$ be an arbitrary bijection of $S$ onto $T$. It is
clear that $(\theta\tsq\theta)\vert_{\cPA(S)}$ is a $\cPA$-isomorphism of $S$ onto $T$ precisely when for all $\alpha\in\cI_S$, we have $\alpha\in\cPA(S)$ if and only if $\alpha(\theta\tsq\theta)\in\cPA(T)$. In particular, any isomorphism or anti-isomorphism of $S$ onto $T$ induces a $\cPA$-isomorphism of $S$ onto $T$.

A semigroup $S$ is called $\cPA${\em-determined} if $S$ is isomorphic or anti-isomorphic to a semigroup $T$ whenever $T$ is
$\cPA$-isomorphic to $S$. We say that $S$ is {\em strongly} $\cPA${\em-determined} if each $\cPA$-isomorphism of $S$ onto a
semigroup $T$ is induced by an isomorphism or an anti-isomorphism of $S$ upon $T$. Let $\K$ be a certain class of semigroups. The
$\cPA${\em -closure} of $\K$ is the class $\cPA(\K)$ of semigroups such that $T\in\cPA(\K)$ if and only if $T$ is
$\cPA$-isomorphic to some $S\in\K$. We say that $\K$ is $\cPA${\em-closed} if $\cPA(\K)=\K$, that is, if $\K$ contains every
semigroup which is $\cPA$-isomorphic to some semigroup from $\K$.

\br{\rm (A
corollary to [19, Main Theorem and its proof].)} Let $S$ be a semilattice (that is, an inverse semigroup such that $S=E_S$) and
$T$ an arbitrary semigroup. Then $\cPA(S)\cong\cPA(T)$ if and only if $S\cong T$ or $S$ is a chain and $T\cong S^d$. Moreover,
any $\cPA$-isomorphism $\Phi$ of $S$ onto $T$ is induced by the $E$-bijection $\varphi_{E}$ associated with $\Phi$, and
$\varphi_{E}$ is either an isomorphism or, if $S$ is a chain and $T\cong S^d$, a dual isomorphism of $S$ onto $T$. \er In
addition, it is easy to see that if a bijection $\gamma$ of $S\,(=E_S)$ onto $T\,(=E_T)$ induces a $\cPA$-isomorphism $\Phi$ of
$S$ onto $T$, then $\gamma=\varphi_{E}$.

Let $S$ be an inverse semigroup. Since we assume that $\emptyset\in\cSbi(S)$, the set of all {\em
inverse} subsemigroups of $S$, partially ordered by inclusion, is a complete lattice which we will denote by $\cSbi(S)$. It is clear that $\cSbi(S)$ is a sublattice of $\cSb(S)$. Among all the partial automorphisms of $S$
it is natural to distinguish those which are isomorphisms between {\em inverse} subsemigroups of $S$; we call them {\em partial
i-automorphisms} of $S$. Since $\emptyset\in\cSbi(S)$, we can also regard $\emptyset$ as a partial {\em i}-automorphism of $S$.
Denote by $\cPAi\,(S)$ the set of all partial {\em i}-automorphisms of $S$. It is clear that $\cPAi\,(S)$ is closed under
composition. Actually, $\cPAi\,(S)$ is an inverse monoid with zero which is an inverse submonoid of $\cPA(S)$. The idempotents
of $\cPAi\,(S)$ are the identity mappings $1_H$ for $H\in\cSbi(S)$, and $E_{\cPAi\,(S)}\cong\cSbi(S)$.

Let $S$ and $T$ be inverse semigroups. If $\cSbi(S)\cong\cSbi(T)$, then $S$ and $T$ are said to be {\em projectively
isomorphic}, and any isomorphism of $\cSbi(S)$ onto $\cSbi(T)$ is called a {\em projectivity} of $S$ upon $T$ (here we use the
terminology of [21]). Again it is clear that an inverse subsemigroup $U$ of $S$ is an atom of $\cSbi(S)$ if and only if $U=\lis
e\ris=\{e\}$ for some $e\in E_S$. Thus if $\Psi$ is a projectivity of $S$ onto $T$, there is a unique  bijection $\psi_E$ of
$E_S$ onto $E_T$ defined by the formula $\{e\}\Psi=\{e\psi_E\}$ for all $e\in E_S$, and we say that $\psi_E$ is the {\em
$E$-bijection associated with $\Psi$}. Since $\cSb(E)=\cSbi(E)$ for any semilattice $E$, there is no difference between lattice
isomorphisms and projectivities of semilattices. In particular, it is again immediate (and well-known) that $\psi_E$ is a weak
isomorphism of $E_S$ onto $E_T$. An important role in the study of projectivities of inverse semigroups is played by the
following result established by Jones in [9]: \br{\rm (From [9, Proposition 1.6 and Corollary 1.7])} Let $S$ and $T$ be
projectively isomorphic inverse semigroups, and let $\Psi$ be a projectivity of $S$ onto $T$. Then there is a (unique) bijection
$\psi : \Ng S\cup E_S\rightarrow \Ng T\cup E_T$ with the following properties:\\ (a) $\psi$ extends $\psi_E$, that is,
$\psi\vert_{E_S}=\psi_E$;\\(b) $\psi$ and $\psi^{-1}$ preserve $\cR$- and $\cL$-classes;\\ (c) $\lis x\ris\Psi=\lis x\psi\ris$
for every $x\in \Ng S\cup E_S$;\\ (d) if a homomorphism $\gamma : S\rightarrow T$ induces $\Psi$, then $x\psi=x\gamma$ for all
$x\in \Ng S\cup E_S$. \er Following [21], we say that the bijection $\psi : \Ng S\cup E_S\rightarrow \Ng T\cup E_T$ in Result
2.2 is the {\em base partial bijection} associated with the projectivity $\Psi$ of $S$ onto $T$. Recall that a semigroup $S$ is
said to be {\em combinatorial} [16, p. 363] if $\cH=1_S$. Clearly a regular semigroup is combinatorial if and only if it has no
nontrivial subgroups. If $S$ in Result 2.2 is combinatorial, then $T$ is combinatorial as well by [8, Corollary 1.3], and since
in this case, $S=\Ng S\cup E_S$ and $T=\Ng T\cup E_T$, the base partial bijection $\psi$ is actually a bijection of $S$ onto
$T$.

Let $S$ and $T$ be inverse semigroups. If $\cPAi\,(S)\cong\cPAi\,(T)$, then $S$ and $T$ are said to be $\cPAi${\em -isomorphic},
and any isomorphism of $\cPAi\,(S)$ onto $\cPAi\,(T)$ is called a $\cPAi${\em -isomorphism} of $S$ onto $T$. Let $\Phi$ be a
$\cPAi$-isomorphism of $S$ onto $T$. Similarly to the case of $\cPA$-isomorphisms of semigroups, for any $H\in\cSbi(S)$ we
define $H\Phi^*$ by the formula $1_H\Phi=1_{H\Phi^*}$, obtaining a projectivity $\Phi^*$ of $S$ onto $T$. As for
$\cPA$-isomorphisms, the $E$-bijection associated with $\Phi^*$ will be denoted by $\varphi_E$ and said to be {\em associated
with $\Phi$}. Also as for $\cPA$-isomorphisms, we say that a bijection $\varphi\!:S\rightarrow T$ {\em induces} $\Phi$ if for
all $\alpha\in\cPAi\,(S)$ and $x,y\in S$, we have $x\alpha=y$ if and only if $(x\varphi)(\alpha\Phi)=y\varphi$. Again it is
obvious that if $\varphi$ is an isomorphism or an anti-isomorphism of $S$ onto $T$, then $\varphi$ induces a $\cPAi$-isomorphism
of $S$ onto $T$. This time, however, if a $\cPAi$-isomorphism $\Phi$ of $S$ onto $T$ is induced by an anti-isomorphism $\varphi$
of $S$ onto $T$, it is also induced by an isomorphism $\iota_S\circ\varphi$ of $S$ onto $T$ where $\iota_S$ is the {\em natural
involution} on $S$ defined by $\iota_S\!:x\mapsto x^{-1}\;(x\in S)$. \br{\rm [4, Lemma 2.3]} Let $S$ and $T$ be
$\cPA$-isomorphic inverse semigroups, and let $\Phi$ be a $\cPA$-isomorphism of $S$ onto $T$. Then the restriction of $\Phi$ to
$\cPAi\,(S)$ is a $\cPAi$-isomorphism of $S$ onto $T$. \er A statement analogous to the following lemma but dealing with
$\cPAi$-isomorphisms of inverse semigroups was proved in [5, Lemma 7]. Actually, both assertions are special cases of the
corresponding general result about $\cPA$-isomorphisms of algebras of any type, the proof of which is entirely similar to that
of [5, Lemma 7]. \bl Let $S$ and $T$ be $\cPA$-isomorphic semigroups, and let $\Phi$ be an isomorphism of $\cPA(S)$ onto
$\cPA(T)$. Then for each $\alpha\in\cPA(S)$, we have ${\rm dom}\,(\alpha\Phi)=({\rm dom}\,\alpha)\Phi^*$ and ${\rm
ran}\,(\alpha\Phi)=({\rm ran}\,\alpha)\Phi^*$, and hence for any subsemigroup $H$ of $S$, the restriction of $\Phi$ to $\cPA(H)$
is a $\cPA$-isomorphism of $H$ onto $H\Phi^*$. \el A {\em null} semigroup is a semigroup $N$ with zero such that $xy=0$ for all
$x, y\in N$. In what follows we will denote by $N_2$ the $2$-element null semigroup $\{0, z\}$. Let $G$ be an arbitrary group
such that $z\not\in G$, and let $e$ be the identity of $G$. It is plain that the mapping $z\mapsto e$ is a partial homomorphism
of $N^*_2$ to $G$; it determines an extension of $G$ by $N_2$ which we will denote by $G^{\langle 1\rangle}$ and call an {\em
extension of $G$ at the identity by $N_2$}. Thus $G^{\langle 1\rangle}=G\cup\{z\}$ is a semigroup with the operation extending
that of $G$ and such that $z^2=e$ and $zx=xz=x$ for all $x\in G$. (Note that $G^{\langle 1\rangle}$ is an {\em inflation} of $G$
[3, \S\,3.2, Exercise 10] with $e$ being replaced by $\{e, z\}$ and all other elements of $G$ left unchanged.) It is easy to
check directly that $\cPA(C_2)\cong\cPA(N_2)$ (and it is obvious that $N_2\cong C^{\langle 1\rangle}_1$), so $C_2$ and
$C^{\langle 1\rangle}_1$ are $\cPA$-isomorphic. Actually, this simple fact is a special case of part (b) of the following
\br{\rm [11, Theorem 1]} Let $S$
be a monogenic semigroup and $T$ an arbitrary semigroup. Then $\cPA(S)\cong\cPA(T)$ if and only if one of the following holds:\\
{\rm(a)} $S\cong T$; {\rm(b)} $S\cong C_{2n}$ and $T\cong C_n^{\langle 1\rangle}$ for an odd $n\geq 1$; {\rm(c)} $S$ and $T$ are
finite monogenic semigroups such that either $\{S,T\}=\{M(2,2),M(3,1)\}$ or $\{S,T\}=\{M(3,6),M(4,3)\}$. \er This result has the
following immediate corollary: \bl Let $S$ and $T$ be $\cPA$-isomorphic semigroups and $\Phi$ a $\cPA$-isomorphism of $S$ onto
$T$. Then for every $x\in\m S$, there is a unique $y\in\m T$ satisfying $\langle x\rangle\Phi^*=\langle y\rangle$, and the
mapping $\varphi : x\mapsto y$ is a bijection of $\m S$ onto $\m T$, extending $\varphi_E$, such that exactly one of the
following holds: {\rm(a)} ${\rm ind}\,x>1$ and ${\rm ind}\,(x\varphi)>1$; {\rm(b)} $\{\langle x\rangle, \langle
x\varphi\rangle\}=\{C_2, N_2\}$; {\rm(c)} $\langle x\rangle\cong C_2\cong\langle x\varphi\rangle$. Moreover, if $S=\m S$, then
$T=\m T$ and $\varphi$ is the unique bijection of $S$ onto $T$ inducing $\Phi$. \el We will say that the mapping $\varphi$,
described in Lemma 2.6, is the {\em $\Phi$-associated bijection of $\m S$ onto $\m T$} (or, if $S=\m S$, the {\em
$\Phi$-associated bijection of $S$ onto $T$}).

Let $S$ and $T$ be $\cPA$-isomorphic combinatorial inverse semigroups, and let $\Phi$ be a $\cPA$-isomorphism of $S$ onto $T$.
Denote, for short, $\Psi=\Phi\vert_{\cPAi(S)}$. By Result 2.3, $\Psi$ is a $\cPAi$-isomorphism of $S$ onto $T$. Let $\varphi$ be
the $\Phi$-associated bijection of $S$ onto $T$, and let $\psi$ be the base bijection of $S$ onto $T$ associated with the
projectivity $\Psi^*$. It is plain that $\varphi_E=\psi_E$. However, it might happen that for some $x, y\in\Ng S$, we have
$x\varphi=x\psi$ but $y\varphi=y^{-1}\psi$; thus, in general, $\varphi\neq\psi$ and $\varphi\neq\iota_S\circ\psi$ (where
$\iota_S$ is the natural involution of $S$). \med
\section{$\cPA$-closed classes of inverse semigroups}
\med Since $C_2$ and $N_2$ are $\cPA$-isomorphic, the class of inverse semigroups is not $\cPA$-closed. In this section we will
show that this ``anomaly'' is, in a sense, the only one: if we remove from the class of all inverse semigroups those having at
least one isolated subgroup which is a direct product of $C_2$ and a periodic group with no elements of even order, we will
obtain a $\cPA$-closed class of inverse semigroup. It is natural to begin our discussion with groups. \br{\rm [12, Lemma]} Let
$G$ be a group and $S$ a semigroup $\cPA$-isomorphic to $G$. Then either $S$ is a group or $S=Q^{\langle 1\rangle}$ where
$Q$ is a periodic subgroup of $S$ with no elements of even order, which is possible only if $G$ is a periodic group with
a unique $2$-element subgroup. \er This lemma was used in [12] in order to prove that a semigroup $S$ is $\cPA$-isomorphic to an {\em
abelian} group $G$ if and only if either $S\cong G$ or $G$ is a periodic abelian group with a unique $2$-element subgroup $C_2$ and
$S\cong(G/C_2)^{\langle 1\rangle}$ [12, Main Theorem]. This gives a complete description of the $\cPA$-closure of the class of
{\em abelian} groups and shows, moreover, that an abelian group $G$ is $\cPA$-determined if it is not a periodic group with a
unique $2$-element subgroup. The latter result does not hold, of course, for nonabelian groups. However, one can strengthen Result 3.1
and obtain a description of the $\cPA$-closure of the class of {\em all} groups.

Let $G$ be an arbitrary group. Denote by $Z(G)$ the center of $G$. If $G$ has elements of order $2$, they are usually called
{\em involutions}. We will say that $G$ is {\em involution-free} if it has no elements of order $2$. Suppose that $G$ has a
unique $2$-element subgroup $A=\{e, a\}$. It is obvious that $A\subseteq Z(G)$, so $A$ is a normal subgroup of $G$. Assume that $G$ splits over $A$, that is, $G$ has a subgroup $P$ (a complement of $A$ in $G$) such that $A\cap P=\{e\}$ and $AP=G$. In this case, since $A\subseteq Z(G)$, it is clear that $G=A\times P$. Thus a group $G$ with a unique $2$-element subgroup $A$ splits over $A$ if and only if $A$ is a direct factor of $G$. Of course, if $G$ is abelian, $A$ is a direct factor of $G$. If $G$ is a finite group, it is immediate from the Burnside normal complement theorem [18, Theorem 7.50] that $G$ contains a
normal complement $P$ of $A$ and hence $G=A\times P$. However, in general, it is not true that a periodic group $G$ with a unique $2$-element subgroup $A$ splits over $A$. Indeed, as follows from [15, Theorem 31.7], there exists a periodic group $G$ with a unique $2$-element subgroup $A$ such that $A$ is not a direct factor of $G$ (the author is grateful to A. Yu. Ol'shanskii for this remark and
reference). At the same time, according to the following lemma, no such periodic group can be $\cPA$-isomorphic to a semigroup
that is not a group.

\bl Let $G$ be a group, $S$ a semigroup which is not a group, and $\Phi$ a $\cPA$-isomorphism of $G$ onto $S$. Then $S$ contains
a periodic involution-free subgroup $Q$ such that $S=Q^{\langle 1\rangle}$. Let $P=Q(\Phi^{-1})^*$. Then $P$ is an
involution-free periodic subgroup of $G$ and $G=C_2\times P$. \el {\bf Proof.} According to Result 3.1, $S=Q^{\langle 1\rangle}$
is an extension of its involution-free periodic subgroup $Q$ at the identity by the $2$-element null semigroup $N_2=\{0, z\}$,
and $A=\{e, z\}(\Phi^{-1})^*$ is a unique $2$-element subgroup of $G$. Since $Q$ is involution-free, by Lemma 2.4 and Result 2.5, $P$ is
also involution-free and thus $A\cap P=\{e\}$. It is plain that $S=\{e, z\}\vee Q$ whence $G=\{e, z\}(\Phi^{-1})^*\vee
Q(\Phi^{-1})^*=A\vee P=AP$. Therefore $G$ splits over $A$. As mentioned above,
this implies that $G=A\times P$.\\

A special case of the next lemma for abelian groups was proved in [12]. The corresponding part of the proof of [12, Main
Theorem] can be easily adjusted to cover our more general situation. For completeness, we include a full proof modifying some
arguments from [12]. \bl Let $P$ and $Q$ be $\cPA$-isomorphic involution-free periodic groups, let $S=Q^{\langle 1\rangle}$ be
an extension of $Q$ at the identity by the $2$-element null semigroup $N_2=\{0, z\}$, and let $G=A\times P$ where $A=\{e,
a\}\cong C_2$. Then $G$ and $S$ are $\cPA$-isomorphic. \el {\bf Proof.} Take an arbitrary $\alpha\in\cPA(G)$. Suppose that ${\rm ran}\,\alpha\not\subseteq P$. Then $ax\in{\rm ran}\,\alpha$ for some $x\in P$. Since $o(x)$ is odd and $o(a)=2$, we have $a=(ax)^{o(x)}\in {\rm ran}\,\alpha$. Hence $a=a\alpha$ since $a$ is the only involution in $G$. Thus ${\rm dom}\,\alpha\not\subseteq P$. For any $X\subseteq G$, denote by $\alpha_{X}$ the restriction of $\alpha$ to $X\cap {\rm dom}\,\alpha$. We have shown that $\alpha_{P}\in\cPA(P)$ for every $\alpha\in\cPA(G)$. It follows also from the above remarks that $\alpha=\alpha_{P}$ if and only if $a\not\in{\rm dom}\,\alpha$, and if $a\in{\rm dom}\,\alpha$, then $\alpha_{(aP)}\not=\emptyset$ and $\alpha=\alpha_{P}\cup\alpha_{(aP)}$. 

Let $\Psi$ be an arbitrary $\cPA$-isomorphism of $P$ onto $Q$. Define $\alpha\Phi=\alpha\Psi$ if $\alpha=\alpha_{P}$, and $\alpha\Phi=\alpha_{P}\Psi\cup\{(z, z)\}$ if $\alpha\not=\alpha_{P}$. Clearly, $\Phi$ is a bijection of $\cPA(G)$ onto $\cPA(S)$; let us show that it is a $\cPA$-isomorphism of $G$ onto $S$. Take any $\alpha,\,\beta\in\cPA(G)$. If $\alpha=\alpha_{P}$ and $\beta=\beta_{P}$, then
$\alpha\circ\beta=(\alpha\circ\beta)_{P}$ and hence 
$(\alpha\circ\beta)\Phi=(\alpha\circ\beta)\Psi=\alpha\Psi\circ\beta\Psi=\alpha\Phi\circ\beta\Phi$. Assume that
$\alpha\not=\alpha_{P}$ and $\beta=\beta_{P}$. Then
$\alpha\circ\beta=(\alpha_{P}\cup\alpha_{(aP)})\circ\beta_{P}=(\alpha_{P}\circ\beta_{P})\cup(\alpha_{(aP)}\circ\beta_{P})
=\alpha_{P}\circ\beta_{P}$
since $aP\cap P=\emptyset$. Thus
$(\alpha\circ\beta)\Phi=(\alpha_{P}\circ\beta_{P})\Phi=(\alpha_{P}\circ\beta_{P})\Psi=\alpha_{P}\Psi\circ\beta_{P}\Psi$, and  \[\alpha\Phi\circ\beta\Phi=(\alpha_{P}\Psi\cup\{(z,
z)\})\circ\beta_{P}\Psi=(\alpha_{P}\Psi\circ\beta_{P}\Psi)\cup(\{(z, z)\}\circ\beta_{P}\Psi)=\alpha_{P}\Psi\circ\beta_{P}\Psi\]
since $z\not\in{\rm dom}\,(\beta_{P}\Psi)\subseteq Q$. Therefore  $(\alpha\circ\beta)\Phi=\alpha\Phi\circ\beta\Phi$. Similarly,
$(\alpha\circ\beta)\Phi=\alpha\Phi\circ\beta\Phi$ if $\alpha=\alpha_{P}$ and $\beta\not=\beta_{P}$. Suppose, finally, that
$\alpha\not=\alpha_{P}$ and $\beta\not=\beta_{P}$. Since $aP\cap P=\emptyset$, we have
\[\alpha\circ\beta=(\alpha_{P}\cup\alpha_{(aP)})\circ(\beta_{P}\cup\beta_{(aP)})=
(\alpha_{P}\circ\beta_{P})\cup(\alpha_{(aP)}\circ\beta_{(aP)}).\]
Hence \[(\alpha\circ\beta)\Phi=(\alpha_{P}\circ\beta_{P})\Psi\cup\{(z, z)\}=(\alpha_{P}\Psi\circ\beta_{P}\Psi)\cup\{(z, z)\}.\]
Using the fact that $z\not\in{\rm dom}\,(\beta_{P}\Psi)$ and $z\not\in{\rm ran}\,(\alpha_{P}\Psi)$, we  also obtain
\[\alpha\Phi\circ\beta\Phi=(\alpha_{P}\Psi\cup\{(z, z)\})\circ(\beta_{P}\Psi\cup\{(z,
z)\})=(\alpha_{P}\Psi\circ\beta_{P}\Psi)\cup\{(z, z)\}.\] Therefore $(\alpha\circ\beta)\Phi=\alpha\Phi\circ\beta\Phi$ in this case as
well. The lemma is proved. \\

A complete description of the $\cPA$-closure of the class of all groups is an immediate consequence of the following  proposition obtained by combining Lemmas 3.2 and 3.3:
 
\bp A group $G$ is $\cPA$-isomorphic to a semigroup $S$ that
is not a group if and only if $G=A\times P$ and $S=Q^{\langle 1\rangle}$ where $P$ and $Q$ are $\cPA$-isomorphic involution-free periodic subgroups of $G$ and $S$, respectively, and $A\cong C_2$. \ep 

The next natural step is to consider $\cPA$-isomorphisms
of monogenic inverse semigroups. \br{\rm (From [13, Main Theorem and its proof]} Let $S=\langle x, x^{-1}\rangle$ be a monogenic
inverse semigroup which is neither a group nor a bicyclic semigroup, and let $T$ be an arbitrary semigroup. Then
$\cPA(S)\cong\cPA(T)$ if and only if $S\cong T$. More precisely, let $\Phi$ be a $\cPA$-isomorphism of $S$ onto $T$, let
$\varphi$ be the $\Phi$-associated bijection of $\m S$ onto $\m T$, and let $y=x\varphi$ and $z=x^{-1}\varphi$. Then $T=\langle
y, z\rangle$ is a monogenic inverse semigroup with $z=y^{-1}$. Moreover, $\varphi$ extends to a bijection $\widetilde{\varphi} :
S\rightarrow T$ which is either an isomorphism of $S$ onto $T$ if $(xx^{-1})\varphi=yy^{-1}$, or an anti-isomorphism of $S$ onto
$T$ if $(xx^{-1})\varphi=y^{-1}y$. \er It should be noted that in [13] the Main Theorem stated that if $A$ is a monogenic
inverse semigroup which is not a group and $B$ is any semigroup, then $\cPA(A)\cong\cPA(B)$ if and only if $A\cong B$. However
the proof of that theorem (see [13], the proof of Lemma 1 on page 55) was based on an erroneous assertion that if $A=\langle a,
a^{-1}\rangle$ is a monogenic inverse semigroup which is not a group, then $A\setminus\{a,a^{-1},aa^{-1},a^{-1}a\}$ is an ideal
of $A$. That assertion is true if, in addition, it is assumed that $A$ is not a bicyclic semigroup. Thus we had to state the
quoted theorem with an additional assumption that the given monogenic inverse semigroup is not bicyclic (that is, how it had
actually been proved in [13]). On the other hand, it is immediate from the proof given in [13], that if $\Phi$ is a
$\cPA$-isomorphism of a monogenic inverse semigroup $S$ onto a semigroup $T$ and $S$ is neither a group nor a bicyclic
semigroup, then there is a bijection $\widetilde{\varphi}$, extending the $\Phi$-associated bijection $\varphi$ of $\m S$ onto
$\m T$, which is either an isomorphism or an anti-isomorphism of $S$ onto $T$, and we stated the quoted theorem in that slightly
sharper form. \br{\rm (From [20, Main Theorem and its proof].)} Let $S=\cB(x, x^{-1})$ be a bicyclic semigroup, $T$ a semigroup
and $\Psi$ a lattice isomorphism of $S$ onto $T$. Then $T$ is also a bicyclic semigroup and $\Psi$ is induced by a bijection
$\psi : S\rightarrow T$, uniquely determined by $\Psi$, such that either $T=\cB(x\psi, (x\psi)^{-1})$, in which case $\psi$ is
an isomorphism, or $T=\cB((x\psi)^{-1}, x\psi)$, in which case $\psi$ is an anti-isomorphism of $S$ onto $T$. \er Now we can
prove the following generalization of Result 3.5: \bl Let $S=\langle x, x^{-1}\rangle$ be a monogenic inverse semigroup, which is
not a group, and $T$ an arbitrary  semigroup. Then $S$ and $T$ are $\cPA$-isomorphic if and only if they are isomorphic.
More specifically, let $\Phi$ be a $\cPA$-isomorphism of $S$ onto $T$, and let $\varphi$ be the $\Phi$-associated bijection of
$\m S$ onto $\m T$. Then there is a bijection $\widetilde{\varphi} : S\rightarrow T$ such that $\widetilde{\varphi}\vert_{\m S}=\varphi$, $T=\langle x\varphi, (x\varphi)^{-1}\rangle$ is a monogenic inverse semigroup, and $\widetilde{\varphi}$ is either an isomorphism of $S$ onto $T$ if $(xx^{-1}) \varphi=(x\varphi)(x\varphi)^{-1}$, or an anti-isomorphism of $S$ onto $T$ if $(xx^{-1})\varphi=(x\varphi)^{-1}(x\varphi)$. Furthermore, if $S=\m S$, then $\varphi$ is the unique bijection of $S$ onto $T$ inducing $\Phi$. \el 

{\bf Proof.} If $S$ is not a bicyclic semigroup, all statements of the lemma (except the one in the last
sentence) follow from Result 3.5, and if $S=\m S$, then $S$ is a ``$C$-semigroup'' in the terminology of [21], so according to [21, Lemma 31.5], the lattice isomorphism $\Phi^*$ is induced by a unique bijection $\varphi$ (which, in this case, coincides
with $\widetilde\varphi$). Now assume that $S=\cB(x, x^{-1})$ is a bicyclic semigroup. Since $\Phi^*$ is a lattice isomorphism
of $S$ onto $T$, by Result 3.6, $T$ is also a bicyclic semigroup and $\Phi^*$ is induced by a unique bijection of $S$ onto $T$
which obviously coincides with $\varphi$. Moreover, either $T=\cB(x\varphi, (x\varphi)^{-1})$ and $\varphi$ is an isomorphism,
or $T=\cB((x\varphi)^{-1}, x\varphi)$ and $\varphi$ is an anti-isomorphism of $S$ onto $T$.

It remains to show that if $S=\m S$, then $\varphi$ induces $\Phi$. Thus suppose that $S=\m S$ (this holds, of course, if $S$ is
combinatorial and, in particular, if $S$ is bicyclic). Take an arbitrary $\alpha\in\cPA(S)$ and any $(x,y)\in\alpha$. Set
$\alpha_x=\alpha\vert_{\langle x\rangle}$. Then $\alpha_x$ is an isomorphism of $\langle x\rangle$ onto $\langle y\rangle$ and
$\alpha_x\subseteq\alpha$. By Lemma 2.4, $\alpha_x\Phi$ is an isomorphism of $\langle x\rangle\Phi^*$ onto $\langle
y\rangle\Phi^*$. Since $\varphi$ induces $\Phi^*$, we have $\langle x\rangle\Phi^*=\langle x\varphi\rangle$ and $\langle
y\rangle\Phi^*=\langle y\varphi\rangle$. Therefore $\alpha_x\Phi$ is an isomorphism of $\langle x\varphi\rangle$ onto $\langle
y\varphi\rangle$ and hence $(x\varphi, y\varphi)\in\alpha_x\Phi\subseteq\alpha\Phi$. Considering $\Phi^{-1}$ and using the
(obvious)  fact that $(\Phi^{-1})^*$ is induced by $\varphi^{-1}$, we obtain, by symmetry, that if
$(x\varphi, y\varphi)\in\alpha\Phi$, then $(x, y)\in\alpha$. This completes the proof.\\

Let $S$ be an inverse semigroup and $T$ an arbitrary semigroup $\cPA$-isomorphic to $S$. Let $\Phi$ be a $\cPA$-isomorphism of
$S$ onto $T$ and $\varphi$ the $\Phi$-associated bijection of $\m S$ onto $\m T$. It is clear that $E_T\not=\emptyset$ and
$\varphi_E\,(=\varphi\vert_{E_S})$ is a bijection of $E_S$ onto $E_T$. By Lemma 2.4, $\Phi\vert_{\cPA(E_S)}$ is a
$\cPA$-isomorphism of $E_S$ onto $E_S\Phi^*(=E_T)$. Hence, according to Result 2.1, $\varphi_E$ is an isomorphism or, if $E_S$
is a chain, perhaps a dual isomorphism of $E_S$ onto $E_T$. In any case, $E_T$ is a semilattice, that is, $T$ is an
idempotent-commutative semigroup. Therefore $\rg (T)$ is the largest inverse subsemigroup of $T$. In what follows, we will
denote $\rg (T)$ by $V$ and $V(\Phi^{-1})^*$ by $U$, so $\Phi\vert_{\cPA(U)}$ is a $\cPA$-isomorphism of $U$ onto $V$. It is
clear that $E_V=E_T$ and $K^V_v=K_v$ for any $\cK\in\{\cH, \cL, \cR, \cD\}$ and $v\in V$. To simplify notation, we will also set
$\Psi=\Phi^{-1}$ (so $\Psi\vert_{\cPA(V)}$ is a $\cPA$-isomorphism of $V$ onto $U$), and denote by $\psi$ the $\Psi$-associated
bijection of $\m T$ onto $\m S$. Let $x$ be an arbitrary nongroup element of $S$. Then $\lis x\ris$ is a monogenic inverse
semigroup which is not a group. Set $\Phi_x=\Phi\vert_{\cPA(\lis x\ris)}$. By Lemma 3.7, there is a bijection $\varphi_x$ of
$\lis x\ris$ onto $\lis x\varphi\ris\,(=\lis x\ris\Phi^*)$ which extends $\varphi\vert_{\m{\lis x\ris}}$. Moreover, according to
Lemma 3.7, $\varphi_x$ is either an isomorphism or an anti-isomorphism of $\lis x\ris$ onto $\lis x\varphi\ris$; it is an
isomorphism if $(xx^{-1})\varphi =(x\varphi)(x\varphi)^{-1}$, and an anti-isomorphism if
$(xx^{-1})\varphi=(x\varphi)^{-1}(x\varphi)$. Since $\lis x\varphi\ris$ is not a group, $D_{x\varphi}$ is a regular nongroup
$\cD$-class of $T$ and hence $D_{x\varphi}=D_{x\varphi}^V$. {\em The notation and observations of this paragraph will be used, frequently without further explanation, throughout the rest of this section.} From the above discussion, using also 
$\Psi\vert_{\cPA(V)}$ instead of $\Phi$, we obtain \bl If $y\in\Ng V$, then $\lis y\psi\ris\,(\cong\lis y\ris)$ is a monogenic
inverse semigroup which is not a group, so $y\psi\in\Ng U$. It follows that $\Ng S\cup E_S=\Ng U\cup E_U$ and $\varphi\vert_{\Ng
U\cup E_U}$ is a bijection of $\Ng U\cup E_U$ onto $\Ng V\cup E_V$. Furthermore, if $e\in E_S$, then $D_e$ is a nongroup
$\cD$-class of $S$ if and only if $D_{e\varphi}$ is a regular nongroup $\cD$-class of $T$ (that is, a nongroup $\cD$-class of
$V$), and for all $g\in E_S$, $g\in D_e$ if and only if $g\varphi\in D_{e\varphi}$; in particular, $e$ is isolated in $S$ if and
only if $e\varphi$ is isolated in $T$.\el We will show later that $H_e\cap U$ is a subgroup of $H_e$ for any $e\in E_S$.
Together with the fact that $\Ng S\cup E_S=\Ng U\cup E_U$, this will imply that $U$ is an inverse subsemigroup of $S$ and hence,
by Result 2.3, $\Phi\vert_{\cPAi(U)}$ is a $\cPAi$-isomorphism of $U$ onto $V$. However, as noted at the end of Section 2,
$\varphi\vert_{\Ng U\cup E_U}$ may still be different from the base partial bijection associated with the projectivity
$(\Phi\vert_{\cPAi(U)})^*$ of $U$ onto $V$. \bl For any idempotent $e$ of $S$, either $H_e\Phi^*=H_{e\varphi}$ or
$H_e\Phi^*=H^{\langle 1\rangle}_{e\varphi}\,=H_{e\varphi}\cup\{z\}$, and in the latter case, $z\not\in\rg (T)$. \el {\bf Proof.}
Let $e\in E_S$ and $f=e\varphi$. By Lemma 2.4 and Result 3.1, either $H_e\Phi^*$ is a subgroup of $T$ or $H_e\Phi^*=Q^{\langle
1\rangle}=Q\cup\{z\}$ where  $Q$ is a subgroup of $T$. Suppose the latter holds. Assume that $z\in \rg (T)$. Since $f=z^2\in\lis
z\ris$, it is immediate that $\lis z\ris$ is not a group (otherwise, $z=zf=f\not=z$, a contradiction). Hence, by Lemma 3.7,
$\lis z\psi\ris\,(=\lis z\ris\Psi^*)$ is a monogenic inverse subsemigroup of $S$ which is not a group. However, it is obvious that
$z\psi\in H_e$ which implies that $\lis z\psi\ris$ is a group. This contradiction shows that $z\not\in\rg (T)$.

It is clear that $f\in Q$, so $Q$ is a subgroup of $H_f$. Applying the above argument to $\Psi\vert_{\cPA(V)}$ and using the fact that $S$ has no nonregular elements, we conclude that $H_f\Psi^*$ is a subgroup of $H_e$. Hence $Q\subseteq
H_f\subseteq H_e\Phi^*=Q\cup\{z\}$. If $Q$ were properly contained in $H_f$, we would have $z\in H_f$, which is impossible since $z\not\in\rg (T)$. Therefore $Q=H_f$ and $H_e\Phi^*=H^{\langle 1\rangle}_f$.

Finally, assume that $H_e\Phi^*$ is a subgroup of $T$. Denote $H_e\Phi^*$ by $K$. Clearly $f\in K$, so $K$ is a subgroup of
$H_f$. As above, we see that $H_f\Psi^*$ is a subgroup of $H_e$ and hence $H_f\subseteq H_e\Phi^*=K\subseteq H_f$. Therefore
$K=H_f$, that is, $H_e\Phi^*=H_f$. This completes the proof. \bl If $e$ is an arbitrary nonisolated idempotent of $S$, then
$H_e\Phi^*=H_{e\varphi}$. Therefore if $H_e\Phi^*=H^{\langle 1\rangle}_{e\varphi}$ for some idempotent $e$ of $S$, then $e$ is
isolated. \el {\bf Proof.} Let $e\in E_S$ be nonisolated. By Lemma 3.9, to prove that $H_e\Phi^*=H_{e\varphi}$, we only need to
show that $H_e\Phi^*\subseteq\rg (T)$. Since $e$ is nonisolated, $D_e$ contains an idempotent $g\not=e$. Let $s$ be an arbitrary
element of $H_e$. Take any $a\in R_e\cap L_g$ and $b\in R_g\cap L_e$. Then $H_a H_b=H_e$ by [3, Theorem 2.2.17], so $s=xy$ for
some $x\in H_a$ and $y\in H_b$. Therefore $s\in\lis x, y\ris=\lis x\ris\vee\lis y\ris$ whence $\lis s\ris\subseteq \lis
x\ris\vee\lis y\ris$. According to [3, Theorem 2.18], $x^{-1}\in H_b$ and $y^{-1}\in H_a$. Hence, by [3, Lemma 2.12],
$xx^{-1}=y^{-1}y=e$ and $x^{-1}x=yy^{-1}=g$, so that $\lis x\ris$ and $\lis y\ris$ are monogenic inverse subsemigroups of $S$
which are not groups. By Lemma 3.7, $\lis x\ris\Phi^*\cong\lis x\ris$ and $\lis y\ris\Phi^*\cong\lis y\ris$. Since $\lis
x\ris\Phi^*$ and $\lis y\ris\Phi^*$ are inverse subsemigroups of $T$, they are contained in $\rg (T)$. Thus \(\lis
s\ris\Phi^*\subseteq(\lis x\ris\vee\lis y\ris)\Phi^*=\lis x\ris\Phi^*\vee\lis y\ris\Phi^*\subseteq\rg (T),\) and hence
\[H_e\Phi^*=(\bigvee_{s\in H_e}\lis s\ris)\Phi^*=\bigvee_{s\in H_e}(\lis s\ris\Phi^*)\subseteq \rg (T).\] The second assertion
of the lemma follows immediately from the first and from Lemma 3.9. \bl Suppose that $T$ is not an inverse semigroup. Let $z$ be
any nonregular element of $T$, and let $a=z\psi$. Then $\langle a\rangle$ is an isolated subgroup of $S$ isomorphic to $C_2$.
Let $e$ denote the identity of $\langle a\rangle$ (that is, $e=a^2$). Then $H_{e\varphi}$ and $H_{e\varphi}\Psi^*$ are
involution-free periodic groups, $H_e\Phi^*=H_{e\varphi}\cup\{z\}$ is an extension of $H_{e\varphi}$ at the identity by the
$2$-element null semigroup $N_2=\{0, z\}$, and $H_e=\langle a\rangle\times (H_{e\varphi}\Psi^*)$. \el {\bf Proof.} Since
$z\not\in\rg(T)$, it is clear that $\text{ind }z>1$ and hence $z\in\m T$. Suppose that $a\in\Ng S$. Then, by Lemma 3.7, $\lis
a\ris\Phi^*\,(\cong\lis a\ris)$ is a monogenic inverse semigroup, so $z$, being an element of $\lis a\ris\Phi^*$, has an
inverse, contradicting the assumption that $z\not\in\rg(T)$. Thus $a$ is a group element of $S$. Assume that $\text{ind
}a=\infty$. Then $\lis a\ris$ is an infinite cyclic group, so $\lis a\ris\Phi^*\cong\lis a\ris$ by [21, Lemma 34.8]. Hence $z$,
as an element of the (infinite cyclic) group $\lis a\ris\Phi^*$, has an inverse, again contradicting the assumption that
$z\not\in\rg(T)$. Therefore $\text{ind }a=1$. Since $\text{ind }z>1$ and $a=z\psi$, by Lemma 2.6, $\langle a\rangle\cong C_2$
and $\langle z\rangle\cong N_2$. Moreover, $H_e\Phi^*\not=H_{e\varphi}$ because $z\in H_e\Phi^*$ and $z\not\in\rg(T)$. Thus, by
Lemma 3.9, $H_e\Phi^*=H^{\langle 1\rangle}_{e\varphi}$ and, by Lemma 3.10, $e$ is an isolated idempotent (and so $\langle
a\rangle$ is an isolated subgroup) of $S$. Since $H^{\langle 1\rangle}_{e\varphi} =H_{e\varphi}\cup Z^*_2$ is an extension of
$H_{e\varphi}$ at the identity by $N_2$ and since the nonzero element of $N_2$ is the only nonregular element of $H^{\langle
1\rangle}_{e\varphi}$, we have $N_2=\{0, z\}$ where $z$ is the given nonregular element of $T$. The remaining assertions of
the lemma follow from Lemma 3.2.\\

We summarize some of the results obtained so far in the following proposition which is an immediate consequence of Lemmas 3.8 --
3.11. \bp Let $e$ be an arbitrary idempotent of $S$. If $e$ is nonisolated, then $D^U_e=D_e$. If $e$ is isolated, then either
$H_e\Phi^*=H_{e\varphi}$, in which case $H_e=H^U_e$, or $H_e\Phi^*=H^{\langle 1\rangle}_{e\varphi}$, in which case
$H_e=A_e\times H^U_e$ where $A_e\cong C_2$ and the group $H^U_e$ is periodic and involution-free. It follows that $U$ is an
inverse subsemigroup of $S$, and if $H_e\Phi^*=H_{e\varphi}$ for all isolated idempotents $e$ of $S$, then $S=U$ and $T=V$ so,
in particular, $T$ is an inverse semigroup. \ep  Now we can establish the main result of this section. \bt Let $S$ be an inverse
semigroup such that no maximal isolated subgroup of $S$ is a direct product of $C_2$ and an involution-free periodic group. Let
$T$ be an arbitrary semigroup $\cPA$-isomorphic to $S$. Then $T$ is also an inverse semigroup in which no maximal isolated
subgroup is a direct product of $C_2$ and an involution-free periodic group. Thus the class of all inverse semigroups, in which
no maximal isolated subgroup is a direct product of a periodic involution-free group and the $2$-element cyclic group, is
$\cPA$-closed. \et {\bf Proof.} Recall that we use the notations fixed in the paragraph preceding Lemma 3.8. In particular,
$\Phi$ denotes a $\cPA$-isomorphism of $S$ onto $T$, and $\Psi=\Phi^{-1}$. Let $e$ be any isolated idempotent of $S$. By
assumption, $H_e$ is not a direct product of $C_2$ and a periodic involution-free group. Hence, according to Proposition 3.12,
$H_e\Phi^*=H_{e\varphi}$. Since $e$ is an arbitrary isolated idempotent of $S$, by Proposition 3.12, $T$ is an inverse
semigroup.

Suppose that $f$ is an isolated idempotent of $T$ such that $H_f=B\times Q$ where $B\cong C_2$ and $Q$ is a periodic
involution-free group. Let $e=f\psi$. In view of Lemma 3.8, $e$ is an isolated idempotent of $S$, and $H_f\Psi^*=H_e$ by Lemma
3.9. Let $A=B\Psi^*$ and $P=Q\Psi^*$. According to Result 2.5, $A\cong C_2$. Since $Q$ is a periodic group, it follows from [1, Theorem 3.2] that $P$ is also periodic, and by Result 2.5, $P$ is involution-free. Since \[ H_e=H_f\Psi^*=(B\vee
Q)\Psi^*=(B\Psi^*)\vee(Q\Psi^*)=A\vee P,\] we have $H_e=A\times P$, which contradicts the condition imposed on $S$. Therefore no maximal isolated subgroup of $T$ is a direct product of $C_2$ and an involution-free periodic group.\\

From this theorem, we can easily deduce that various classes of inverse semigroups are $\cPA$-closed (in the class of all
semigroups). For example, we have
\bc The following classes of inverse semigroups are $\cPA$-closed:\\
(a) the class of all inverse semigroups with no isolated subgroups of order $2$;\\
(b) the class of all inverse semigroups with no nontrivial isolated subgroups;\\
(c) the class of all combinatorial inverse semigroups. \ec {\bf Proof.} Recall again that we use the notation of the paragraph
preceding Lemma 3.8. Suppose that $S$ has no isolated subgroups of order $2$. Then no maximal isolated subgroup of $S$ is a
direct product of $C_2$ and a periodic involution-free group. Thus, by Theorem 3.13, $T$ is an inverse semigroup. Assume that
$T$ contains an isolated subgroup $B$ of order $2$ and denote by $f$ the identity of $B$. Let $e=f\psi$ and $A=B\Psi^*$. Then
$e$ is an isolated idempotent of $S$, and it is clear that $A\cong C_2$, so $A$ is an isolated subgroup of $S$ of order $2$; a
contradiction. This proves (a), whereas (b) and (c) follow immediately from Theorem 3.13 and Lemma 3.9. \\

We conclude this section with an example of a class of Clifford semigroups which are $\cPA$-isomorphic to semigroups that are
not inverse. Recall that a {\em Clifford semigroup} is a regular semigroup in which the idempotents are central. The structure
of Clifford semigroups was completely determined in [2] by means of the following construction. Let $E$ be an arbitrary
semilattice, and let $S_e\;(e\in E)$ be a family of pairwise  disjoint semigroups. Suppose that for all $e,f\in E$ with $e\geq
f$, there is a homomorphism $\varphi_{e,f}\!:\,S_e\rightarrow S_f$ such that $\varphi_{e,e}=1_{S_e}$ and
$\varphi_{e,f}\circ\varphi_{f,g}=\varphi_{e,g}$ for all $e,f,g\in E$ satisfying $e\geq f\geq g$. If we define multiplication on
$S=\bigcup\{S_e\,\vert\,e\in E\}$ by the formula $s\ast t=(s\varphi_{e,ef})(t\varphi_{f,ef})$ for all $s, t\in S$ (where $s\in
S_e,\,t\in S_f$), then $(S,\ast)$ becomes a semigroup called a {\em strong semilattice} $E$ {\em of semigroups} $S_e$ {\em
determined by the homomorphisms} $\varphi_{e,f}$, which is written as $S=[E;S_e,\varphi_{e,f}]$ (see [16, II.2.2 and II.2.3]).
In [2] Clifford proved that $S$ is a regular semigroup with central idempotents if and only if $S$ is a strong semilattice of
groups.

Let $E$ be an arbitrary semilattice, and let $A=[E;A_e,\varphi_{e,f}]$ where $A_e=\{e,a_e\}\cong C_2$ for each $e\in E$ and
$A_e\varphi_{e,f}=\{f\}$ for all $e,f\in E$ such that $e > f$. Let $B=[E;B_e,\psi_{e,f}]$ where $B_e=\{e,z_e\}\cong N_2$ for
each $e\in E$ and $B_e\psi_{e,f}=\{f\}$ for all $e,f\in E$ such that $e > f$. Now let $\theta$ be a bijection of $A$ onto $B$
such that $e\theta=e$ and $a_e\theta=z_e$ for every $e\in E$. Let $\Theta=(\theta\tsq\theta)\vert_{\cPA(A)}$. It is easily seen
that if $\alpha\in\cI_A$, then $\alpha\in \cPA(A)$ if and only if $\alpha(\theta\tsq\theta)\in\cPA(B)$, and thus $\Theta$ is an
isomorphism of $\cPA(A)$ onto $\cPA(B)$. In short, we have \bex Let $A,\,B,\,\theta$, and $\Theta$ be as defined
in the preceding paragraph. Then $A$ is a Clifford semigroup, $B$ is a combinatorial semigroup which is not inverse, and
$\Theta$ is a $\cPA$-isomorphism of $A$ onto $B$ induced by $\theta$. \eex Thus the class of Clifford semigroups (which are not groups) is not $\cPA$-closed and neither is the class of (nontrivial) combinatorial semigroups. Using Example 3.15
(and its modifications) and taking Lemma 3.11 as a starting point, we can obtain a complete description of the $\cPA$-closure of the class of all inverse semigroups, which will be given in another article.\med
\section{$\cPA$-determined inverse semigroups}
\med In this section we consider the problem of $\cPA$-determinability of inverse semigroups. Let $S$ be an inverse semigroup. If
$a\in S$ and $e\in E_S$ are such that $e<aa^{-1}$ and there is no $f\in E_{[\![a]\!]}$ satisfying  $e<f<aa^{-1}$, we say that
$e$ is $a${\em -covered} by $aa^{-1}$. Take any $a\in S$ and $e\in E_S$ with $e<aa^{-1}$. Suppose that for some positive integer
$n$, there exist $e_{0}, e_{1}, \ldots, e_n\in E_S$ such that $e=e_{0}<e_{1}<\cdots<e_n=aa^{-1}$ and for every $k=1,\ldots,n$,
the idempotent $e_{k-1}$ is $a_{k}$-covered by $e_{k}$ where $a_{k}=e_{k}a$ (and hence $a_{k} a^{-1}_{k}=e_{k}$). Then
$(e_{0},e_{1},\ldots,e_n)$ is called a {\em short bypass} from $e$ to $aa^{-1}$. If for all $e,\,a\in S$ such that $e<aa^{-1}$,
there is a short bypass from $e$ to $aa^{-1}$, then $S$ is said to be a {\em shortly connected} inverse semigroup. This property
was introduced in [5] in connection with the following theorem:

\br {\rm [5, Theorem 5]} Let $S$ be a combinatorial inverse semigroup, $T$ an inverse semigroup projectively isomorphic to $S$,
and $\Psi$ a projectivity of $S$ onto $T$. Let $\psi$ be the base bijection of $S$ onto $T$ associated with $\Psi$ (so, in
particular, $\psi_{E}=\psi\vert_{E_S}$). Suppose that $S$ is shortly connected and $\psi_{E}$ is an isomorphism of $E_S$ onto
$E_T$. Then $\psi$ is the unique isomorphism of $S$ onto $T$ which induces $\Psi$. \er An inverse semigroup $S$ is called {\em
shortly linked} if for all $a\!\in\!S$ and $e\!\in\!E_S$ such that $e<aa^{-1}$, the set $F_{e,a}=\{f\in E_{[\![a]\!]}: e<f\leq
aa^{-1}\}$ is finite. By [5, Proposition 3], any shortly linked inverse semigroup is shortly connected. In fact, the class of
shortly linked inverse semigroups is properly contained in the class of shortly connected ones [6]. However, the property of
being shortly linked is easier to check than the property of being shortly connected, and precisely for this reason shortly
linked inverse semigroups were introduced in [5]. Thus it might be useful to formulate the following specialization of Result
4.1 (as an obvious consequence of [5, Theorem 5], it was not explicitly stated in [5]): \br {\rm (A corollary to [5, Theorem
5])} Let $S$ be a combinatorial inverse semigroup, $T$ an inverse semigroup projectively isomorphic to $S$, and $\Psi$ a
projectivity of $S$ onto $T$. Let $\psi$ be the base bijection of $S$ onto $T$ associated with $\Psi$. Suppose $S$ is
shortly linked and $\psi_{E}$ is an isomorphism of $E_S$ onto $E_T$. Then $\psi$ is the unique isomorphism of $S$ onto $T$inducing $\Psi$. \er

In general, if $\Psi$ is a projectivity of a combinatorial inverse semigroup $S$ onto an inverse semigroup $T$, the
$E$-bijection $\psi_E$ associated with $\Psi$ may not be an isomorphism but just a weak isomorphism of $E_S$ onto $E_T$;
however, in many interesting special cases $\psi_E$ is, in fact, an isomorphism of $E_S$ onto $E_T$ (see [5] and [10] for more
details). The original reason for imposing this condition in Results 4.1 and 4.2 was the fact that if $\Phi$ is a
$\cPAi$-isomorphism of a combinatorial inverse semigroup $S$ onto an inverse semigroup $T$, then the base bijection $\varphi$ of
$S$ onto $T$ associated with $\Phi^*$ is such that $\varphi_E$ is indeed an isomorphism of $E_S$ onto $E_T$, except for the case
when $(S, \leq)$ is a chain, $T=S^d$, and $\varphi_E\,(=\varphi)$ is a dual isomorphism of $S$ onto $T$. More precisely, using
Result 4.1, we proved in [5] the following theorem:

\br {\rm [5, Theorem 8]} Let $S$ be a shortly connected combinatorial inverse semigroup and $T$ an inverse semigroup. Then
$\cPAi(S)\cong\cPAi(T)$ if and only if either $S\cong T$ or $(S,\leq)$ and $(T,\leq)$ are dually isomorphic chains. Moreover,
any $\cPAi$-isomorphism of $S$ onto $T$ is induced by a unique isomorphism of $S$ onto $T$ or, if $(S,\leq)$ is a chain and
$T\cong S^d$, by a unique dual isomorphism of $S$ onto $T$. \er

Let $S$ be an inverse semigroup. As indicated above, the requirement that $S$ be shortly connected is strictly weaker than the
requirement that it be shortly linked. Several other properties of $S$ which are strictly weaker than the property of being
shortly linked were introduced recently in [10]. Following Jones [10], we will call $S$ {\em pseudo-archimedean} if none of its
idempotents is strictly below every idempotent of a bicyclic or free monogenic inverse subsemigroup of $S$, {\em faintly
archimedean} if whenever an idempotent $e$ of $S$ is strictly below every idempotent of a bicyclic or free monogenic inverse
subsemigroup $\lis a\ris$ of $S$ then $e < a$, and {\em quasi-archimedean} if it is faintly archimedean and $\lis x\ris$ is
combinatorial for each $x\in\Ng S$. Thus it is immediate that every pseudo-archimedean inverse semigroup is faintly archimedean.
Note that Jones defined the property of $S$ being quasi-archimedean differently and then proved it to be equivalent to the one
given above in [10, Proposition 3.3(3)], which makes apparent the fact that for combinatorial inverse semigroups the properties
of being faintly archimedean and quasi-archimedean coincide [10, Corollary 3.4].) According to [10, Example 3.2], there exists a
pseudo-archimedean inverse semigroup $S$ which is also combinatorial and {\em $E$-unitary} (that is, such that for any $a\in S$
and $e\in E_S$, if $e\leq a$, then $a\in E_S$) but which is not shortly linked. For quasi-archimedean combinatorial inverse
semigroups, the following theorem was proved in [10]: \br {\rm [10, Theorem 4.3]} Let $S$ be a combinatorial inverse semigroup,
$T$ an inverse semigroup projectively isomorphic to $S$, and $\Psi$ a projectivity of $S$ onto $T$. Let $\psi$ be the base
bijection of $S$ onto $T$ associated with $\Psi$. Suppose that $S$ is quasi-archimedean (equivalently, faintly archimedean) and
$\psi_{E}$ is an isomorphism of $E_S$ onto $E_T$. Then $\psi$ is the unique isomorphism of $S$ onto $T$ which induces $\Psi$.
\er As pointed out in [10], this theorem generalizes Result 4.2 since every shortly linked inverse semigroup is faintly
archimedean. However, the question of whether every shortly connected inverse semigroup is faintly archimedean was not addressed
in [10]. Now we will construct two examples of shortly connected combinatorial inverse semigroups which are not faintly
archimedean (one of them will contain a bicyclic subsemigroup while the other one will be completely semisimple), showing
therefore that Result 4.4 does not generalize Result 4.1.

Recall that an inverse semigroup $S$ is said to be {\em fundamental} if $1_S$ is the only congruence on $S$ contained in $\cH$,
so every combinatorial inverse semigroup is certainly fundamental. Fundamental inverse semigroups, introduced by Munn [14] (and
independently by Wagner [24] under a different name), form one of the most important classes of inverse semigroups (see [7] and
[14] for details). Let $E$ be an arbitrary semilattice. Recall that the {\em Munn semigroup} $T_E$ is an inverse semigroup
(under composition) consisting of all isomorphisms between principal ideals of $E$ [7, \S\,V.4]. If $S$ is an inverse semigroup,
a subset $K$ of $S$ is called {\em full} if $E_S\subseteq K$ [16, p. 118]. Munn proved (see [14, Theorem 2.6] or [7, Theorem
V.4.10]) that an inverse semigroup $S$ with $E_S=E$ is fundamental if and only if $S$ is isomorphic to a full inverse
subsemigroup of $T_E$, and hence $T_E$ itself is fundamental.

Let \( E=\{e_0, e_1, e_2, \ldots, f_0, f_1, f_2,\ldots, g_0, g_1, 0\}\) be the semilattice given by the diagram in Figure 1. Let
$S=T_E$ be the Munn semigroup of the semilattice $E$. As usual, we will identify each $e\in E$ with $1_{Ee}\in E_S$ (so that
$E_S$ is identified with $E$). It is immediate that $Ee_m\cong Ee_n$ and $Ef_m\cong Ef_n$ for all integers $m, n\geq 0$, and
$Ee_m\ncong Ef_m$ and $Ee_m\ncong Eg_0\cong Eg_1\ncong Ef_m$ for every $m\geq 0$. In fact, it is easy to see that for any
integers $m, n\geq 0$, there is exactly one isomorphism $\varphi_{m, n}$ of $Ee_m$ onto $Ee_n$; it is defined as follows:
\(e_k\varphi_{m, n}=e_{k-m+n}\;\text{ and }\;f_k\varphi_{m, n}=f_{k-m+n}\;\text{ for all integers }\;k\geq m,\) $g_0\varphi_{m,
n}=g_0$ and $g_1\varphi_{m, n}=g_1$ if $m-n\equiv 0\,({\rm mod\;}2)$, $g_0\varphi_{m, n}=g_1$ and $g_1\varphi_{m, n}=g_0$ if
$m-n\equiv 1\,({\rm mod\;}2)$, and $0\varphi_{m, n}=0$. It is clear that $(\varphi_{m, n})^{-1}=\varphi_{n, m}$ and the
restriction of $\varphi_{m, n}$ to $Ef_m$ is the only isomorphism of $Ef_m$ onto $Ef_n$ (while the restriction of $\varphi_{m,
n}$ to $Eg_0$ or to $Eg_1$ is also the only isomorphism between the corresponding principal ideals of $E$).
\vspace{0.05in}\\
\begin{picture}(360,237)
\put(228,54){\line(-1,-1){34}} \put(160,52){\line(1,-1){33}}
\put(40,100){\line(5,4){155}}\put(342,99){\line(-5,-2){113}}\put(303,99){\line(-5,-3){75}}\put(258,100){\line(-2,-3){30}}
\put(80,100){\line(5,4){115}}\put(120,100){\line(5,4){75}}\put(40,99){\line(5,-2){118}}\put(80,100){\line(5,-3){78}}
\put(120,100){\line(5,-6){39}} \put(340,100){\line(-4,3){143}}\put(300,100){\line(-4,3){103}}\put(256,100){\line(-4,3){61}}
\put(195,224){\line(0,-1){14}}\put(195,208){\line(0,-1){14}}\put(195,192){\line(0,-1){14}}\put(195,178){\line(0,-1){16}}
\put(195,160){\line(0,-1){14}}\put(195,145){\line(0,-1){10}} \put(195,224){\circle*{5}}\put(195,208){\circle*{5}}
\put(195,192){\circle*{5}} \put(195,178){\circle*{5}}\put(195,160){\circle*{5}}\put(195,145){\circle*{5}}
\put(193,19){\circle*{5}} \put(227,53){\circle*{5}} \put(159,53){\circle*{5}}\put(40,100){\circle*{5}}
\put(80,100){\circle*{5}}\put(120,100){\circle*{5}}\put(340,100){\circle*{5}}
\put(301,100){\circle*{5}}\put(257,100){\circle*{5}}\put(199,10){$0$} \put(200,225){$e_{_0}$} \put(183,203){$e_{_1}$}
\put(200,190){$e_{_2}$}\put(183,172){$e_{_3}$}\put(200,156){$e_{_4}$}\put(183,138){$e_{_5}$}\put(232,47){$g_{_1}$}
\put(145,47){$g_{_0}$}\put(25,100){$f_{_0}$}\put(65,100){$f_{_2}$}\put(105,100){$f_{_4}$}\put(345,100){$f_{_1}$}
\put(305,100){$f_{_3}$}\put(261,100){$f_{_5}$} \put(140,100){$\ldots$}\put(220,100){$\ldots$}\put(195,110){$\vdots$}
\put(175,-10){$\text{Figure }1$}
\end{picture}
\vspace{0.1in}\\

As observed in [7, the proof of Proposition V.6.1], if $F$ is an arbitrary semilattice, then for any $e, f\in F$, we have $(e,
f)\in \cD$ in $T_F$ if and only if $Fe\cong Ff$. It follows that $S$ is a combinatorial inverse semigroup with exactly four
$\cD$-classes: $D_0$, $D_{g_{_0}}$, $D_{f_{_0}}$, and $D_{e_{_0}}$; moreover, $E_{D_0}=D_0=\{0\}$, $E_{D_{g_{_0}}}=\{g_0,
g_1\}$, $E_{D_{f_{_0}}}=\{f_0, f_1,\ldots\}$, and $E_{D_{e_{_0}}}=\{e_0, e_1, \ldots\}$. Let $a=\varphi_{0, 1}$. It is obvious
that $D_{e_{_0}}$ is the bicyclic semigroup $\cB(a, a^{-1})$ and $D_0\cup D_{g_{_0}}\cup D_{f_{_0}}$ is a completely semisimple
inverse subsemigroup of $S$ containing the five-element Brandt subsemigroup $D_0\cup D_{g_{_0}}$. Note that $g_0 < e_m$ for all
$m\geq 0$ (that is, each idempotent of $\cB(a, a^{-1}$) is strictly above $g_0$) but $g_0\nless a$. This means that $S$ is not
faintly archimedean. On the other hand, it is easily seen that $S$ is shortly connected. Thus we have the following \bex Let $E$
be the semilattice whose diagram is shown in Figure 1. Then the Munn semigroup $T_E$ is a shortly connected combinatorial
inverse semigroup which contains a bicyclic subsemigroup and is not faintly archimedean. \eex

Now let $E$ be a semilattice whose diagram is shown in Figure $2$. Its subsemilattice \(E'=\{e_{10}, e_{01}; e_{20}, e_{11},
e_{02}; e_{30}, e_{21}, e_{12}, e_{03}; \ldots\}\) is the semilattice of idempotents of the free monogenic inverse semigroup
where, as in [6], $\{e_{n-q,\,q} : q=0, 1, \ldots, n\}\;(n\in\N)$ is the set of all idempotents of that semigroup of weight $n$
[16, Sections IX.1 and IX.2]. Again as in [6], $e_{pq}$ stands here for the idempotent that can be uniquely written in the form
$e_pf_q$ in the notation of [16, p. 408] where $p, q\geq 0$ and $p+q>0$. Furthermore, our semilattice $E$ contains a primitive
subsemilattice $\{g_0, g_1, 0\}$ and pairwise incomparable elements $f_{pq}$ (where $p+q=n$ and $n$ runs through the set of all
odd positive integers) such that every $f_{pq}$ is covered in $E$ by $e_{pq}$, and $f_{pq}$, in turn, covers either $g_0$ if $p$
is odd, or $g_1$ if $p$ is even.
\vspace{0.05in}\\
\begin{picture}(380,237)
\put(210,225){\line(-1,-3){27}}\put(191,224){\line(1,-3){27}}\put(190,225){\line(-3,-2){154}}\put(210,225){\line(3,-2){154}}
\put(210,225){\line(1,-3){29}}\put(190,225){\line(-1,-3){29}}\put(180,196){\line(1,-3){17}}\put(220,196){\line(-1,-3){17}}
 \put(170,167){\line(1,-3){8}}
\put(230,167){\line(-1,-3){8}} \put(147,80){\line(4,-3){55}}  \put(254,80){\line(-4,-3){55}} \put(83,123){\circle*{5}}
\put(37,123){\circle*{5}} \put(363,123){\circle*{5}}\put(317,123){\circle*{5}}\put(278,123){\circle*{5}}
\put(147,80){\circle*{5}}\put(254,80){\circle*{5}}\put(37,123){\line(5,-2){110}}\put(83,123){\line(3,-2){64}}
\put(122,123){\line(3,-5){26}} \put(363,123){\line(-5,-2){110}}\put(317,123){\line(-3,-2){64}}
\put(279,123){\line(-3,-5){26}}\put(170,167){\line(-2,-1){87}} \put(230,167){\line(2,-1){87}}\put(210,167){\line(-2,-1){87}}
\put(190,167){\line(2,-1){87}} \put(190,225){\circle*{5}}\put(210,225){\circle*{5}}\put(200,196){\circle*{5}}
\put(122,123){\circle*{5}}\put(200,40){\circle*{5}}
\put(180,196){\circle*{5}}\put(220,196){\circle*{5}}\put(170,167){\circle*{5}}\put(190,167){\circle*{5}}
\put(210,167){\circle*{5}}\put(230,167){\circle*{5}} \put(205,30){$0$} \put(172,227){$e_{_{10}}$} \put(215,227){$e_{_{01}}$}
\put(163,197){$e_{_{20}}$}\put(203,197){$e_{_{11}}$}\put(225,197){$e_{_{02}}$}\put(153,169){$e_{_{30}}$}\put(175,169){$e_{_{21}}$}
\put(213,169){$e_{_{12}}$} \put(232,169){$e_{_{03}}$} \put(20,119){$f_{_{10}}$}\put(66,119){$f_{_{30}}$}
\put(104,119){$f_{_{12}}$} \put(284,119){$f_{_{21}}$}\put(320,119){$f_{_{03}}$} \put(366,119){$f_{_{01}}$}\put(240,80){$g_{_1}$}
\put(153,80){$g_{_0}$} \put(159,120){$\ldots$}\put(230,120){$\ldots$}\put(200,130){$\vdots$} \put(182,7){$\text{Figure }2$}
\end{picture}
\vspace{0.01in}

Let $\alpha$ be the isomorphism of $Ee_{10}$ onto $Ee_{01}$ that is uniquely determined by the formula $e_{pq}\alpha=e_{p-1,
q+1}$ ($p\geq 1, q\geq 0$). Note that whenever $p\geq 1$ and $q\geq 0$ are such that $p+q$ is odd, then $f_{pq}\alpha=f_{p-1,
q+1}$. Furthermore, $g_0\alpha=g_1$, $g_1\alpha=g_0$, and $0\alpha=0$. Let $S$ be the full inverse subsemigroup of $T_E$
generated by $\alpha$. It is easily seen that $S$ is shortly connected and combinatorial. Since $\lis \alpha\vert_{E'}\ris$ is
the free monogenic inverse subsemigroup of $S$ and since $g_0 < \alpha\alpha^{-1}\,(=e_{10})$ but $g_0\nless \alpha$, we
conclude that $S$ is not faintly archimedean. Note that $S$ does not contain a bicyclic subsemigroup, so it is completely
semisimple. Thus we have \bex Let $E$ be the semilattice whose diagram is shown in Figure 2. Let $\alpha$
be an isomorphism of $Ee_{10}$ onto $Ee_{01}$ uniquely determined by the formula $e_{pq}\alpha=e_{p-1, q+1}$ for any $p\geq 1,
q\geq 0$, and let $S$ be the full inverse subsemigroup of the Munn semigroup $T_E$ generated by $\alpha$. Then $S$ is a
completely semisimple shortly connected combinatorial inverse semigroup which is not faintly archimedean. \eex

For any semigroup $S$, denote by $S^0$ the semigroup obtained from $S$ by adjoining an ``extra'' zero element $0$ to $S$. It is clear that if $S$ is either a free
monogenic inverse semigroup or a bicyclic semigroup, then $S^0$ is faintly archimedean but not shortly connected. Together with the above two examples, this shows that the properties of being shortly connected and faintly
archimedean for inverse semigroups are independent of one another.

 \bt Let $S$ be a combinatorial inverse semigroup which is either shortly connected or faintly archimedean, and let $T$ be an
arbitrary semigroup. Then $\cPA(T)\cong\cPA(S)$ if and only if $T$ is an inverse semigroup such that either $T\cong S$ or
$(S,\leq)$ is a chain and $(T,\leq)\cong(S,\leq^{d})$. More specifically, if $\Phi$ is a $\cPA$-isomorphism of $S$ onto $T$, it
is induced by a unique bijection $\varphi$ of $S$ onto $T$ such that either $(S,\leq)$ and $(T, \leq)$ are dually isomorphic
chains and $\varphi\,(=\varphi_E)$ is a dual isomorphism of $S$ onto $T$, or $\varphi_E$ is an isomorphism of $E_S$ onto $E_T$,
in which case $\varphi\vert_{\lis x\ris}$ is an isomorphism or an anti-isomorphism of $\lis x\ris$ onto $\lis x\varphi\ris$ for
every $x\in \Ng S$, and the $\cPAi$-isomorphism $\Phi\vert_{\cPAi(S)}$ of $S$ onto $T$ is induced by a unique isomorphism of $S$
onto $T$. \et {\bf Proof.} Let $\Phi$ be a $\cPA$-isomorphism of $S$ onto $T$. According to Corollary 3.14(c), $T$ is a
combinatorial inverse semigroup. Let $\varphi$ be the $\Phi$-associated bijection of $S$ onto $T$. Then, in particular,
$\varphi\vert_{E_S}=\varphi_E$ where $\varphi_E$ is the $E$-bijection associated with $\Phi$. By Result 2.1, either $\varphi_E$
is an isomorphism of $E_S$ onto $E_T$, or $(E_S, \leq)$ and $E_T, \leq)$ are dually isomorphic chains and $\varphi_E$ is an
isomorphism of $E_S$ onto $(E_T)^d$. Suppose that the latter holds. Then, since $S$ and $T$ are combinatorial, it is not
difficult to show that $S=E_S$ and $T=E_T$ (see the last paragraph of the proof of [5, Theorem 8]), so that
$\varphi\,(=\varphi_E)$ is the unique bijection of $S$ onto $T$ inducing $\Phi$, and $\varphi$ is a dual isomorphism of $(S,
\leq)$ onto $(T, \leq)$. Now assume that $\varphi_E$ is an isomorphism of $E_S$ onto $E_T$. From Result 2.3 and Lemma 3.7, it
follows that $\varphi$ is the unique bijection of $S$ onto $T$ inducing $\Phi$, and $\varphi\vert_{\lis x\ris}$ is an
isomorphism or an anti-isomorphism of $\lis x\ris$ onto $\lis x\varphi\ris$ for every $x\in\Ng S$. Finally, by Result 2.3,
$\Phi\vert_{\cPAi(S)}$ is a $\cPAi$-isomorphism of $S$ onto $T$, and hence $(\Phi\vert_{\cPAi(S)})^*$ is a projectivity of $S$
onto $T$. Therefore, if $S$ is either shortly connected or quasi-archimedean, then $S\cong T$ by Result 4.1 or by Result 4.4,
respectively. Using Result 2.3, we obtain the last assertion of the theorem by applying Result 4.3 in case $S$ is
shortly connected, and deduce it from Result 4.4 if $S$ is faintly archimedean.\\

Under the assumptions and in the notation of Theorem 4.7, it is natural to ask: Is it true that the $\Phi$-associated bijection
$\varphi$ of $S$ onto $T$ (in the case when $\varphi_E$ is an isomorphism of $E_S$ onto $E_T$) is either an isomorphism or an
anti-isomorphism of $S$ onto $T$? In general, the answer is no. For example, let $A=\langle a, a^{-1}\rangle$ be the free
monogenic inverse semigroup, let $B=\{0, b, b^{-1}, bb^{-1}, b^{-1}b\}$ be the five-element Brandt semigroup, and let $S$ be an extension of $B$ by $A^0$ determined by the map $a\mapsto b$. Then $S$ is a faintly archimedean combinatorial inverse semigroup. Define a bijection $\varphi : S\rightarrow S$ as follows: $b\varphi=b^{-1}$, $b^{-1}\varphi=b$, and $\varphi\vert_{S\setminus\{b, b^{-1}\}}=1_{S\setminus\{b, b^{-1}\}}$, and let $\Phi=(\varphi\tsq\varphi)\vert_{\cPA(S)}$. Note that for every $s\in S\setminus E_S$, if $s\in A$ then $\lis s\ris$ is a free monogenic inverse subsemigroup of $S$, and if $s\in B$ then $\lis s\ris=B$. It follows that for an arbitrary $\alpha\in\cI(S)$, we have $\alpha\in\cPA(S)$ if and only if $\alpha\Phi\in\cPA(S)$. Therefore $\Phi$ is a $\cPA$-isomorphism of $S$ onto $T$. By the very definition, $\Phi$ is induced by $\varphi$. However, $\varphi$ is neither an isomorphism nor an anti-isomorphism of $S$ onto $S$ because $\varphi\vert_A$ is an isomorphism of $A$ onto $A$ whereas $\varphi\vert_B$ is an anti-isomorphism of $B$ onto
$B$. Using the same idea, one can construct other faintly archimedean inverse semigroups with analogous properties. Thus we have the following \bex There exist faintly archimedean (combinatorial) inverse semigroups $S$ such that there is a
$\cPA$-isomorphism of $S$ onto an inverse semigroup $T\,(\cong S)$ induced by a (unique) bijection which is neither an
isomorphism nor an anti-isomorphism of $S$ onto $T$. \eex

There are, of course, similar examples of shortly connected inverse semigroups. Nevertheless for some classes of combinatorial inverse semigroups $S$ every $\cPA$-isomorphism of $S$ onto a semigroup $T$ is induced either by an isomorphism or by an anti-isomorphism of $S$ onto $T$. As an illustration, we will give one example of such class. 
\bp Let $S$ be a periodic combinatorial inverse semigroup, $T$ an arbitrary semigroup $\cPA$-isomorphic to $S$, and  $\Phi$ any $\cPA$-isomorphism of $S$ onto $T$. Then $\Phi$ is induced by a unique bijection $\varphi$ which is either an isomorphism or an anti-isomorphism of $S$ onto $T$. \ep {\bf Proof.} Since $S$ is periodic, it is immediate that $S$ is shortly linked, so it is both shortly connected and faintly archimedean. According to Theorem 4.7, $T$ is an inverse semigroup isomorphic to $S$ and the $\Phi$-associated bijection $\varphi : S\rightarrow T$ has the following properties: $\varphi_E$ is an isomorphism of $E_S$ onto $E_T$ and for every $s\in S\setminus E_S$, the restriction of $\varphi$ to $\lis s\ris$ is an isomorphism or an anti-isomorphism of $\lis s\ris$ onto $\lis s\varphi\ris$. Suppose that $x, y\in\Ng S$ are such that $\varphi\vert_{\lis x\ris}$ is an isomorphism of $\lis x\ris$ onto $\lis x\varphi\ris$ but $\varphi\vert_{\lis y\ris}$ is an anti-isomorphism of $\lis y\ris$ onto $\lis y\varphi\ris$. Since $S$ is periodic, $\lis x\ris$ and $\lis y\ris$ are finite combinatorial inverse subsemigroups of $S$. Denote by $e$ the least idempotent of $\lis x\ris$ and by $f$ the least idempotent of $\lis y\ris$. Take any $a\in\lis x\ris$ and $b\in\lis y\ris$ such that $aa^{-1}\succ e$ and $bb^{-1}\succ f$. Then $\lis a\ris$ is a five-element Brandt subsemigroup of $\lis x\ris$ and $\lis b\ris$ a five-element Brandt subsemigroup of $\lis y\ris$. Moreover, $\varphi\vert_{\lis a\ris}$ is an isomorphism of $\lis a\ris$ onto $\lis a\varphi\ris$, and $\varphi\vert_{\lis b\ris}$ is an anti-isomorphism of $\lis b\ris$ onto $\lis b\varphi\ris$. Let $A=\{a, aa^{-1}, e\}$ and $B=\{b, bb^{-1}, f\}$. It is clear that $\alpha : A\rightarrow B$ given by:
$a\alpha =b$, $(aa^{-1})\alpha=bb^{-1}$, and $e\alpha=f$, is the unique isomorphism of $A$ onto $B$, so $\alpha\in\cPA(S)$.
Hence $\alpha\Phi\in\cPA(T)$, that is, $\alpha\Phi$ is an isomorphism of $A\Phi^*$ onto $B\Phi^*$. At the same time,
$A\Phi^*=\{a\varphi, (a\varphi)(a\varphi)^{-1}, e\varphi\}$ and $B\Phi^*=\{b\varphi, (b\varphi)^{-1}(b\varphi), f\varphi\}$
because $\varphi\vert_A$ is an isomorphism of $A$ onto $A\Phi^*$ whereas $\varphi\vert_B$ is an anti-isomorphism of $B$ onto
$B\Phi^*$. Thus we have
\[b\varphi=(a\varphi)\alpha\Phi=[(a\varphi)(a\varphi)^{-1}\cdot a\varphi]\alpha\Phi=[(a\varphi)(a\varphi)^{-1}]\alpha\Phi
\cdot(a\varphi)\alpha\Phi=(b\varphi)^{-1}(b\varphi)\cdot b\varphi=f\varphi;\] a contradiction. Therefore either
$\varphi\vert_{\lis s\ris}$ is an isomorphism of $\lis s\ris$ onto $\lis s\varphi\ris$ for all $s\in\Ng S$ or
$\varphi\vert_{\lis s\ris}$ is an anti-isomorphism of $\lis s\ris$ onto $\lis s\varphi\ris$ for all $s\in\Ng S$. From this it easily follows that $\varphi$ is either an isomorphism or an anti-isomorphism of $S$ onto $T$.\\

\begin{center}{\bf Acknowledgement}\end{center}
The work on this paper was completed while the author held the position of a Visiting Professor at the Department of Mathematics at Stanford University from January through June of 2004, and he thanks the department's faculty and staff for their hospitality and for creating excellent working conditions. The author is especially grateful to Professor Yakov Eliashberg for the invitation to Stanford and for his hospitality.

\end{document}